\documentclass[12pt]{amsart}
\usepackage{amsfonts,latexsym,amsthm,amssymb,graphicx}
\usepackage[all]{xy}

\DeclareFontFamily{OT1}{rsfs}{}
\DeclareFontShape{OT1}{rsfs}{n}{it}{<-> rsfs10}{}
\DeclareMathAlphabet{\mathscr}{OT1}{rsfs}{n}{it}

\newtheorem{theorem}{Theorem}[section]
\newtheorem{lemma}[theorem]{Lemma}
\newtheorem{corol}[theorem]{Corollary}
\newtheorem{prop}[theorem]{Proposition}

\newtheorem{conj}{Conjecture}

{\theoremstyle{definition} \newtheorem{defin}[theorem]{Definition}}
{\theoremstyle{remark} 
\newtheorem{example}[theorem]{Example}}

\newcommand{\Abb}{{\mathbb{A}}}
\newcommand{\Cbb}{{\mathbb{C}}}

\newcommand{\Pbb}{{\mathbb{P}}}

\newcommand{\Sbb}{{\mathbb{S}}}
\newcommand{\Vbb}{{\mathbb{V}}}
\newcommand{\Zbb}{{\mathbb{Z}}}

\newcommand{\one}{1\hskip-3.5pt1}

\newcommand{\cE}{{\mathscr{E}}}
\newcommand{\cF}{{\mathscr{F}}}

\newcommand{\cL}{{\mathscr{L}}}
\newcommand{\cO}{{\mathscr{O}}}
\newcommand{\cQ}{{\mathscr{Q}}}
\newcommand{\cS}{{\mathscr{S}}}

\DeclareMathOperator{\rk}{rk}

\newcommand{\csm}{c_{\text{\rm SM}}}
\newcommand{\ualpha}{{\underline\alpha}}
\newcommand{\ubeta}{{\underline\beta}}

\title[Chern classes of Schubert cells and varieties]
{Chern classes of Schubert cells and varieties}

\author{Paolo Aluffi}
\author{Leonardo Constantin Mihalcea}
\address{
Mathematics Department, 
Florida State University,
Tallahassee FL 32306, U.S.A.
}
\email{aluffi@math.fsu.edu}
\address{
Mathematics Department, 
Florida State University,
Tallahassee FL 32306, U.S.A.
}
\email{mihalcea@math.fsu.edu}

\begin{document}

\begin{abstract}
We give explicit formulas for the Chern-Schwartz-MacPherson classes of
all Schubert varieties in the Grassmannian of $d$-planes in a vector
space, and conjecture that these classes are effective. We prove this
is the case for (very) small values of $d$.
\end{abstract}

\maketitle


\section{Introduction}\label{intro}

The classical {\em Schubert varieties\/} in the Grassmannian $G_d(V)$
parametrize subspaces of dimension $d$ (`$d$-planes') of an ambient
vector space $V$, satisfying prescribed incidence conditions with a
flag of subspaces. Schubert varieties are among the most studied
objects in algebraic geometry, and it may come as a surprise that
something is {\em not\/} known about them. Yet, to our knowledge the
{\em Chern classes\/} of Schubert varieties have not been available in
the literature. This paper is devoted to their computation.

Most Schubert varieties are singular; in fact, the only
nonsingular Schubert varieties are themselves (isomorphic to)
Grassmannians, and computing their Chern classes is a standard
exercise. For singular varieties, there is a good notion
of Chern classes in `homology' (i.e., in the Chow group), that
is, the so-called {\em Chern-Schwartz-MacPherson\/} (`CSM') classes.
They were defined independently by Marie-H\'el\`ene Schwartz (\cite{MR35:3707},
\cite{MR32:1727}) and Robert MacPherson (\cite{MR50:13587}), they
agree with the total homology Chern class of the tangent bundle for
nonsingular varieties, and they satisfy good functoriality properties
(summarized, for example, in \cite{MR85k:14004}, Example~19.1.7).

It is a consequence of these good functoriality properties
classes that a Chern-Schwartz-MacPherson class can in fact be
assigned to all Schubert {\em cells:\/} these are varieties isomorphic
to affine spaces, parametrizing $d$-planes which satisfy prescribed
incidence conditions {\em strictly.\/} Every Schubert {\em variety\/}
is stratified by the Schubert {\em cells\/} contained in it, and as a
consequence the CSM class of a Schubert variety may be written as a
sum of CSM classes of Schubert cells contained in it.

With this understood, our task amounts to the computation of the CSM
class of a Schubert cell, as an element of the Chow group of a
Schubert variety containing it.

In order to state the result more precisely, we have to bring in some
notation. We consider the `abstract' Schubert variety $\Sbb(\ualpha)$
determined by a partition $\ualpha=(\alpha_1\ge \alpha_2\ge\dots \ge
\alpha_d \ge 0)$; we depict $\ualpha$ by the corresponding
(upside-down) Young diagram~$\ualpha$. For example, $\Sbb(7\ge 5\ge
3\ge 2)$ is the Schubert variety associated to the diagram
\begin{center}
\includegraphics[scale=.6]{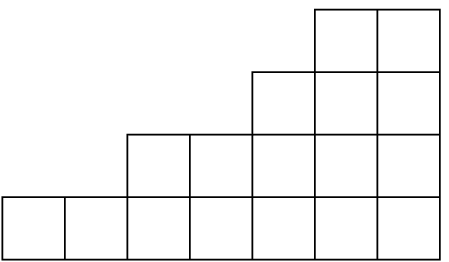}
\end{center}
Embedding this Schubert variety in (for example) the Grassmannian
$G_5(V)$ of $5$-planes in a 13-dimensional vector space $V$ realizes
it as the subvariety parametrizing subspaces intersecting a fixed flag
of subspaces of dimensions~$1$, $4$, $6$, $9$, $12$ in dimension $\ge
1,2,3,4, 5$ respectively. It is perhaps more common to associate this
latter realization with the `complementary' Young diagram,
\begin{center}
\includegraphics[scale=.6]{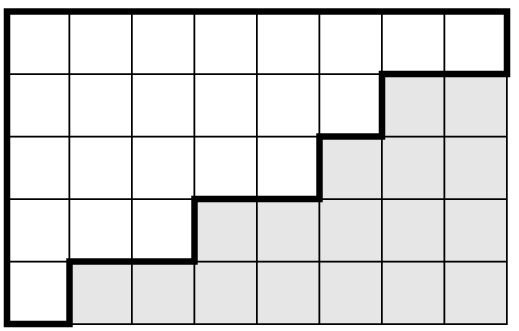}
\end{center}
whose Schur function determines the class of $\Sbb(\ualpha)$ in the 
{\em cohomology\/} (or Chow ring) of $G_5(V)$. However, since CSM classes 
live in {\em homology,\/} we are naturally led to take the dual viewpoint. 
This also has the advantage of allowing us to deal with Schubert 
varieties as abstract varieties, independently of a (standard) 
embedding into a Grassmannian. The dimension of $\Sbb(\ualpha)$ 
equals the number of boxes in the corresponding `homological' diagram.

We write $\ubeta\le\ualpha$ to denote $\beta_i\le\alpha_i$
$\forall i$; that is, the diagram corresponding to $\ubeta$ is
contained in the diagram corresponding to $\ualpha$. For $\ubeta\le
\ualpha$ there are closed embeddings $\Sbb(\ubeta)\subset
\Sbb(\ualpha)$; the Schubert {\em cell\/} $\Sbb(\ualpha)^\circ$ is the
complement of $\cup_{\beta < \alpha} \Sbb(\ubeta)$ in
$\Sbb(\alpha)$. As mentioned above, Schubert cells are isomorphic to 
affine spaces, and there is a decomposition

$$\Sbb(\ualpha)=\coprod_{\beta\le\alpha} \Sbb(\ubeta)^\circ\quad.$$
It also follows that the Chow group $A_*\Sbb(\ualpha)$ is freely 
generated by the classes $[\Sbb(\ubeta)]$ for all $\ubeta\le\ualpha$. 
For $\ubeta\le\ualpha$, the CSM class 
$$\csm(\Sbb(\ubeta)^\circ)\in A_*\Sbb(\ualpha)$$
of the corresponding Schubert cell is defined
as the image via MacPherson's natural correspondence 
(cf.~\cite{MR50:13587}) of the constructible function 
$\one_{\Sbb(\ubeta)^\circ}$ whose value is~1 on $\Sbb(\ubeta)^\circ$,
and $0$ on its complement in $\Sbb(\ualpha)$.
Since $\one_{\Sbb(\ualpha)}=
\sum_{\beta\le \alpha} \one_{\Sbb(\ubeta)^\circ}$,
the basic covariance property of CSM classes implies that
$$\csm(\Sbb(\ualpha)
=\sum_{\ubeta\le \ualpha} \csm(\Sbb(\ubeta)^\circ)\quad,$$
reducing the problem of computing CSM classes of Schubert
{\em varieties\/} to that of computing CSM classes of Schubert
{\em cells.\/} Now,
$$\csm(\Sbb(\ualpha)^\circ)=\sum_{\ubeta\le \ualpha}
\gamma_{\ualpha,\ubeta} [\Sbb(\ubeta)]$$
for uniquely determined coefficients $\gamma_{\ualpha,\ubeta}\in \Zbb$.
Our main result consists of the computation of these coefficients,
and can be stated in several different forms. For example:

\begin{theorem}\label{main}
Let $\ualpha=(\alpha_1\ge \cdots \ge \alpha_d)$, $\ubeta=(\beta_1\ge 
\cdots \ge \beta_d)$. The integer $\gamma_{\ualpha,\ubeta}$ equals
the coefficient of 
$$t_1^{\alpha_1}\cdots t_d^{\alpha_d}\cdot 
u_1^{\beta_1}\cdots u_d^{\beta_d}$$
in the expansion of the rational function
$$\Gamma_d(\underline t,\underline u)=
\frac 1{(t_1^d\cdots t_d^1)(u_1^d\cdots u_d^1)}\cdot
\prod_{1\le i<j\le d} \frac{(t_i-t_j)(u_i-u_j)}{1-2t_j+t_i t_j}\cdot
\prod_{1\le i, j\le d}
\frac{1-t_i}{1-t_i(1+u_j)}$$
as a Laurent polynomial in $\Zbb[[\underline t,\underline u]]$.
\end{theorem}

(Even) more explicitly, the coefficient $\gamma_{\ualpha,\ubeta}$
equals 
$$\sum \det
\left[
\begin{pmatrix}
\alpha_i-\ell_{i+1}^i-\cdots -\ell_d^i \\
\beta_j+(i-j)+\ell_i^1+\dots+\ell_i^{i-1}-\ell_{i+1}^i-\cdots -\ell_d^i 
\end{pmatrix}
\right]_{1\le i,j\le d}$$
where the summation is over the $\binom d2$ integers $\ell_i^k$, $1\le k
<i\le d$, subject to the conditions
$$0\le \ell_{k+1}^k + \cdots + \ell_d^k \le\alpha_{k+1}\quad.$$

Other expressions for the coefficients $\gamma_{\ualpha,\ubeta}$ are
given in~\S\ref{CSM}. If $\beta_i=0$ for $i\ge 2$, that is, $\ubeta=
(r)$ corresponds to a one-row diagram, then $\gamma_{\ualpha,\ubeta}$ equals the coefficient of 
$u^r$~in
$$\prod_{i\ge 1} (1+iu)^{\alpha_i-\alpha_{i+1}}$$
(see Corollary~\ref{onerow} below). 
For 2-row diagrams $\ualpha$, a formula
for $\gamma_{\ualpha,\ubeta}$ in terms of non-intersecting 
lattice paths is given in Theorem~\ref{latticepaths}.

\begin{example}\label{g25}
The Grassmannian $G_2(\Cbb^5)$ parametrizing lines in $\Pbb^4$
is isomorphic to $\Sbb(\ualpha)$ for a $3\times 2$ rectangle $\ualpha$. 
We list here the CSM classes of its Schubert cells, 
denoting by a darkened diagram $\ubeta$ the corresponding class
$[\Sbb(\ubeta)]\in A_*\Sbb(\ualpha)=A_*G_2(\Cbb^5)$, and by a dot the
`empty' diagram (corresponding to the class of a point):
\begin{center}
\includegraphics[scale=.35]{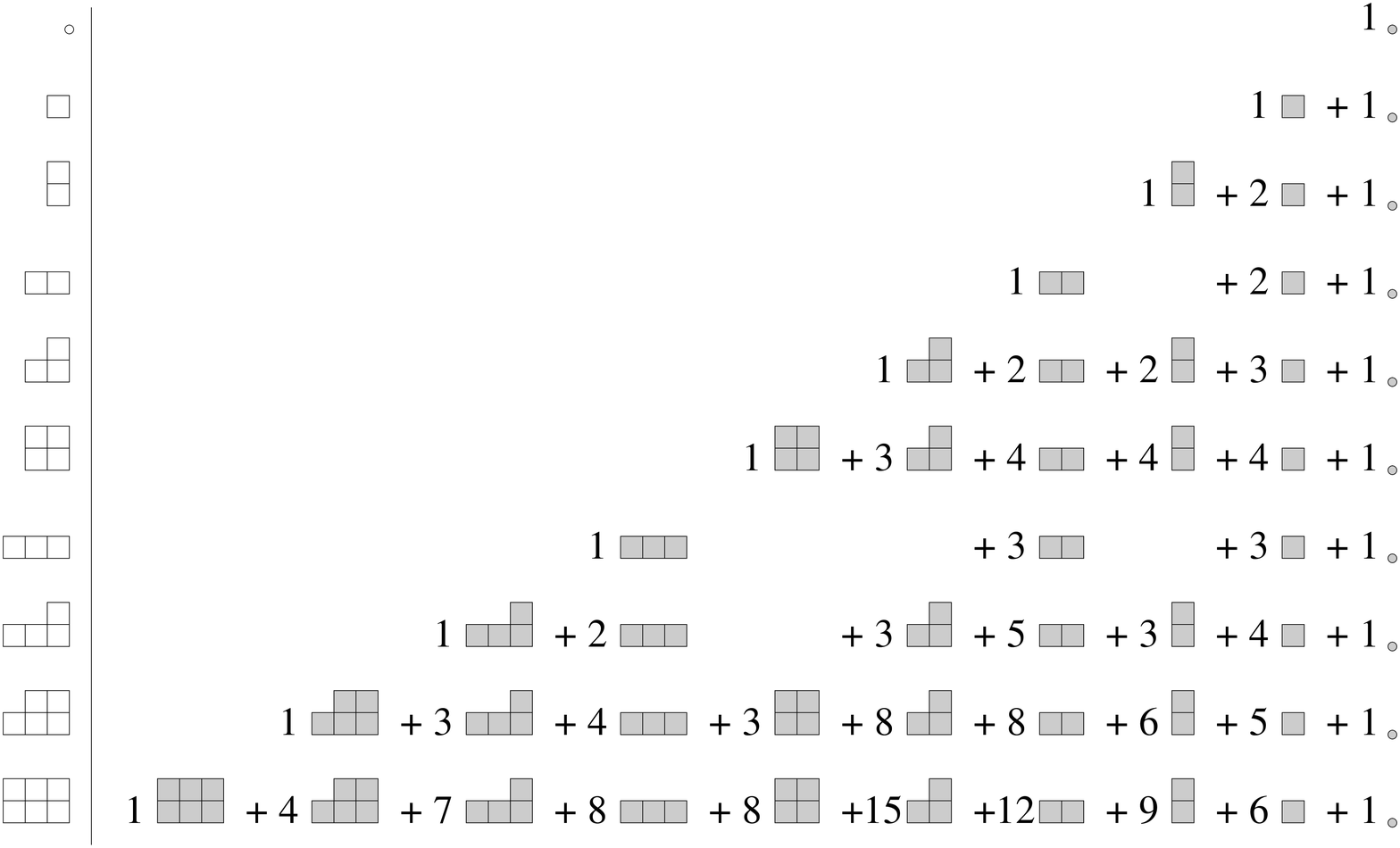}
\end{center}
For example, the Chern-Schwartz-MacPherson class of $G_2(\Cbb^5)$
may be obtained by adding up all these classes:
\begin{center}
\includegraphics[scale=.35]{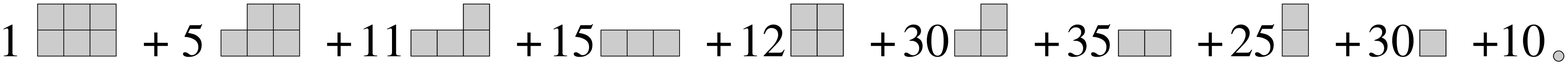}
\end{center}
As $G_2(\Cbb^5)$ is nonsingular, $\csm(G_2(\Cbb^5))$ agrees with the 
total Chern class of its tangent bundle, which may be computed from the 
standard tautological sequence on the Grassmannian.
\end{example}

A number of features of CSM classes of Schubert cells are clear for 
`geometric' reasons. For example:
\begin{itemize}
\item The contribution of the class of a point to the CSM class of a
Schubert cell is~$1$; this is necessarily the case, since each cell is
an affine space and the degree of the CSM class of a variety always
equals its Euler characteristic (by functoriality);
\item {\em Mutatis mutandis,\/} the CSM class for a diagram $\ualpha$
must agree with the CSM class for the transposed diagram, obtained from 
$\ualpha$ by interchanging rows and columns;
\item The CSM class for a sequence $\ualpha=
(\alpha_1\ge\alpha_2\ge\cdots)$
only depends on the nonzero terms in the sequence.
\end{itemize}
These facts correspond to identities involving the coefficients
$\gamma_{\ualpha,\ubeta}$ computed in Theorem~\ref{main}, which may
appear more mysterious from the `algebraic' point of view.

Less immediate features include adjunction-type formulas,
and interpretations of the coefficients $\gamma_{\ualpha,\ubeta}$ as
computing certain nonintersecting lattice paths in the plane; these
are discussed in \S\ref{pos}, and follow from the algebraic
expressions obtained in~\S\ref{CSM}.

There is one feature which is experimentally manifest, and for which
we do not know a general proof from either the algebraic or the
geometric viewpoint: all coefficients of the CSM classes of all
Schubert cells appear to be {\em positive.\/} We give two proofs of
this fact for diagrams with $\le 2$ rows in \S\ref{pos}, as
consequences of the facts mentioned in the previous paragraph.  A
proof of positivity for diagrams with $3$ rows will appear
elsewhere. Positivity of Chern classes is a well-explored theme in
intersection theory; the case of Schubert cells hints that there may
be an interesting, and as yet unknown, principle of positivity for
Chern-Schwartz-MacPherson classes of singular or noncomplete
varieties under suitable hypotheses.

The proof of Theorem~\ref{main} relies on explicit nonsingular
birational models $\Vbb(\ualpha)$ of Schubert varieties
$\Sbb(\ualpha)$---the well-known {\em Bott-Samelson resolutions\/} 
(cf.~\cite{MR0354697}, \cite{math.AG/0302294}).
In \S\ref{BS} we give a self-contained description of these varieties, 
independent of the embedding of $\Sbb(\ualpha)$ in a Grassmannian.
For the application to CSM classes it is necessary to compute
explicitly the push-forward at the level of Chow groups
$$A_*(\Vbb(\ualpha)) \to A_*(\Sbb(\ualpha))\quad;$$
this is accomplished in Proposition~\ref{pfprop},
and may be of independent interest. The varieties $\Vbb(\ualpha)$
realize the splitting principle for the tautological bundle $\cS$ on
$\Sbb(\ualpha)$ (obtained as pull-back of the tautological subbundle
from a Grassmannian), so that the Chern roots of $\cS$ are realized as
Chern classes $-\xi_i$ of line bundles $\cL_i$ on $\Vbb(\alpha)$. 
We explicitly compute the push-forward to $\Sbb(\ualpha)$ of a monomial
in the $\xi_i$'s, as the class (up to sign) of a smaller Schubert variety
$\Sbb(\ubeta)$. The reader will check that the classical Pieri's formula 
is an immediate consequence of this observation.

The variety $\Vbb(\ualpha)$ is a nonsingular compactification of the
Schubert cell $\Sbb(\ualpha)^\circ$, with complement a divisor with
simple normal crossing. The CSM class of $\Sbb(\ualpha)^\circ$ may
then be obtained by computing the (ordinary) Chern class of a bundle
of logarithmic tangent fields on $\Vbb(\ualpha)$. This is done in
\S\ref{CSM}, leading to explicit expressions for
$\csm(\Sbb(\ualpha)^\circ)$, of which Theorem~\ref{main} is a sample.

This approach to the computation of Chern-Schwartz-MacPherson classes
goes back to an observation (Proposition~15.3 in
\cite{MR1893006}, Theorem~1 in \cite{MR1717120}) and has proven
useful for other explicit computations. For example, it gives a
one-sentence proof of Fritz Ehlers' formula for the CSM class of toric
varieties (\cite{MR2209219}).

Chern-Schwartz-MacPherson classes of (possibly) singular varieties
have been the object of intense study, but they have been computed 
explicitly in relatively few cases. To our knowledge, the largest 
class of varieties for which CSM classes are known consists of 
{\em degeneracy loci\/} of morphisms of vector bundles, satisfying
a mild generality hypotheses, which are treated by 
Adam Parusi\'nski and Piotr Pragacz in \cite{MR1311826}. In fact,
there is some overlap of the results in \cite{MR1311826} and in
the present paper: Schubert varieties indexed by 
$\ualpha=(\alpha_1=\alpha_2=\cdots=\alpha_{d-1}\ge \alpha_d)$
(that is, by a `cohomological' one-row diagram) may
be realized as degeneracy loci. Comparing the formulas obtained 
here with those in {\em loc.~cit.\/} is a natural project.

Chern-Schwartz-MacPherson classes are arguably the most natural
replacement for $c(TX)\cap [X]$ when the tangent bundle is not
available, and computing them for varieties as well-known as Schubert
varieties is a natural task. It would also be a natural task to
compute other classes satisfying the same basic normalization
requirement, such as Mather's Chern class (\cite{MR85k:14004},
Example~4.2.9) or Fulton's Chern class (\cite{MR85k:14004},
Example~4.2.6). Comparisons between these different classes for
Schubert varieties would lead to the computation of other
important invariants of their singularities, such as the local 
Euler obstruction.

{\bf Acknowledgements.} We are grateful to the Max-Planck-Institut
f\"ur Mathematik, Bonn, for hospitality and support in the Summer of
2006. We thank Don Zagier for many enlightening discussions on
generating functions, and especially for a series of hints which
catalyzed the proof of Corollary~\ref{onerow}.


\section{Schubert and Bott-Samelson varieties}\label{BS}
\subsection{Schubert varieties.}\label{notat}
We work over an algebraically closed field $k$. In this section we
recall some material concerning classical Schubert varieties, mainly
for the purpose of setting notations. Proofs may be found in any
standard reference, such as \cite{MR85k:14004}, Chapter~14.

We denote by $\ualpha$ a {\em partition,\/} that is, a nonincreasing 
sequence of nonnegative integers $\alpha_i$, such that $\alpha_i=0$ for
$i\gg 0$:
$$\ualpha: (\ualpha_1\ge \ualpha_2\ge\dots\ge 0\ge \dots)\quad.$$
If we write $\ualpha=(\alpha_1\ge\dots\ge \alpha_d)$, it is understood
that $\alpha_i=0$ for $i>d$.

Pictorially, every partition is associated to a Young diagram
\begin{center}
\includegraphics[scale=.6]{pictures/Young}
\end{center}
with $\alpha_i$ boxes in the $i$-th row, counting from the bottom;
with this in mind, we sometime call $\ualpha$ a {\em diagram.\/}

We can associate to $\ualpha$ a {\em Schubert variety\/}
$\Sbb(\ualpha)$, as follows. For $\ualpha=(0)$ (the `empty diagram'),
$\Sbb(\ualpha)$ is a point. For nonzero $\ualpha$, let $N$, $d$ be
integers such that $N\ge \alpha_1$ and $\alpha_{d+1}=0$. Let $V$ be a
vector space of dimension~$N+d$, and fix a complete flag $F_\bullet$
in~$V$:
$$F_0=\{0\}\subset F_1\subset \cdots \subset F_{N+d}=V\quad,$$
where $\dim F_r=r$. The sequence $\ualpha$ selects $d$
elements in the flag:
$$F_{\alpha_d+1}\subset F_{\alpha_{d-1}+2}\subset \cdots
\subset F_{\alpha_1+d}\quad;$$
the dimensions of the spaces in this subflag
are the lengths of the rows of the diagram obtained adjoining to
$\ualpha$ a `$d$-step ladder':
\begin{center}
\includegraphics[scale=.6]{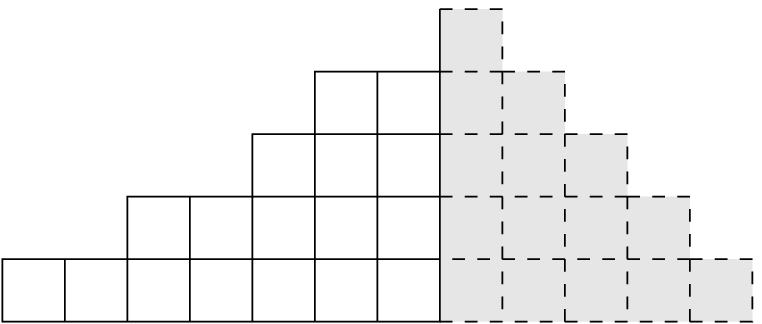}
\end{center}

\begin{defin}\label{defschvar}
The {\em Schubert variety\/} $\Sbb(\ualpha)$ denotes the variety of 
$d$-planes $S$ in the Grassmannian $G_d(V)$ such that 
$$\dim (S\cap F_{\alpha_{d+1-i}+i})\ge i\quad.$$
\end{defin}

Note that $\Sbb(\ualpha)$ is not defined unless $\ualpha$ is a 
nonincreasing sequence of integers, as above.

\begin{lemma}
Up to isomorphism, $\Sbb(\ualpha)$ only depends on $\ualpha$.
\end{lemma}

Indeed, the definition of $\Sbb(\ualpha)$ makes sense as soon as 
$N=\alpha_1$ and $d$ is the number of nonzero elements in $\ualpha$; 
increasing $N$ does not affect the constraints, while increasing $d$ amounts 
to direct-summing a fixed subspace to all subspaces parametrized by the 
variety, producing an isomorphic variety $\Sbb(\ualpha)$.

It is also the case that $\Sbb(\ualpha)$ and $\Sbb(\ualpha^t)$ are
isomorphic, for $\alpha$, $\alpha^t$ `transposed' diagrams,
obtained from each other by interchanging rows and columns.
\begin{center}
\includegraphics[scale=.6]{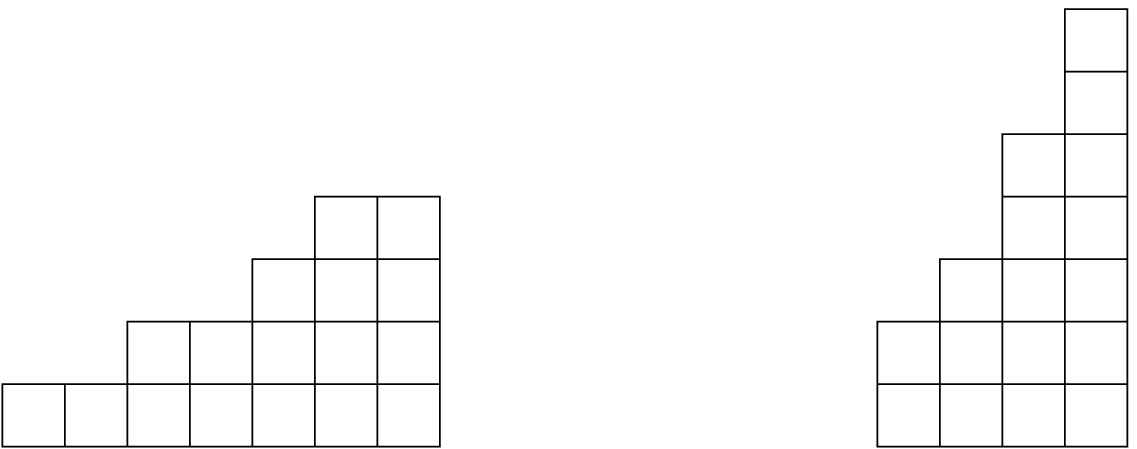}
\end{center}
This is a
consequence of the correspondence between subspaces of dual
vector spaces.

\subsection{Schubert cells.}
We assume that a sufficiently large vector space $V$, a complete
flag $F_\bullet$ in $V$, and a sufficiently large $d$ have been chosen.

We denote by $\ubeta\le\ualpha$ the ordering $\beta_i\le \alpha_i$ 
$\forall i$. If $\ubeta\le \ualpha$, there are embeddings
$$\Sbb(\ubeta)\subset \Sbb(\ualpha)\quad:$$
indeed, the conditions defining $\Sbb(\ubeta)$ are a strengthening 
of those defining $\ualpha$.

\begin{defin}\label{defschcells}
The Schubert {\em cell\/} corresponding to $\ualpha$ is the 
open dense subset
$$\Sbb(\ualpha)^\circ = \Sbb(\ualpha)\setminus
\cup_{\beta < \alpha} \Sbb(\beta)$$
of the Schubert {\em variety\/} $\Sbb(\ualpha)$.
\end{defin}

Equivalently, a $d$-plane $S\in \Sbb(\ualpha)$ is in the cell
$\Sbb(\ualpha)^\circ$ if and only if
$$\dim (S\cap F_{\alpha_{d+1-i}+i})= i\quad\text{and}
\quad \dim(S\cap F_{\alpha_{d+1-i}+i-1})= i-1$$
for $i=1,\dots,d$.

An elementary coordinate computation proves the following:

\begin{lemma}
The Schubert cell $\Sbb(\alpha)^\circ$ is isomorphic to an affine space
$\Abb^{\sum_i \alpha_i}$ of dimension $\sum_i \alpha_i$.
\end{lemma}

It follows that $\dim\,\Sbb(\ualpha)=\sum_i \alpha_i$; further, this
shows that the Schubert varieties have a cellular decomposition
(cf.~\cite{MR85k:14004}, Example~1.9.1) and in particular that
$A_*(\Sbb(\ualpha))$ is freely generated by the {\em Schubert
classes\/} $[\Sbb(\beta)]$ for $\beta\le \alpha$. As stated in the
introduction, our task is to express the Chern-Schwartz-MacPherson
class of the Schubert {\em cell\/} $\Sbb(\ualpha)^\circ$, as a
combination of Schubert classes:
$$\csm(\Sbb(\ualpha)^\circ)=\sum_{\beta\le \alpha}
\gamma_{\ualpha,\ubeta}\, [\Sbb(\ubeta)]\quad \in A_*(\Sbb(\ualpha)]
\quad.$$
 
\subsection{Bott-Samelson varieties}\label{BSvars}
In the embedded setting, there is a natural way to produce a
nonsingular model of the Schubert variety $\Sbb(\ualpha)\subset
G_d(V)$. If $\Sbb(\ualpha)$ is defined with respect to the complete
flag
$$F_0=\{0\}\subset F_1\subset \cdots \subset F_{N+d}=V\quad,$$
as in \S\ref{notat}, let $\Vbb(\ualpha)$ be the corresponding
{\em Bott-Samelson variety\/} of flags
$$\Vbb(\ualpha)=\{ (S^1 \subset S^2 \subset \cdots \subset S^d)\,|\,  
\dim S^i=i \textrm{ and } S^i \subset F_{\alpha_{d+1-i}+i}\}\quad.$$
The top space in the flag, $S^d$, satisfies the conditions defining 
$\Sbb(\ualpha)$; therefore, there is a
natural map
$$\pi_{\alpha}\,:\, \Vbb(\ualpha) \to \Sbb(\ualpha)\quad.$$
It is clear from this description that $\pi_\ualpha$ has a section 
over the Schubert cell $\Sbb(\ualpha)^\circ$, defined by sending
sending $S\in \Sbb(\ualpha)^\circ$ to the flag
$$(S\cap F_{\alpha_{d}+1}\subset S\cap F_{\alpha_{d-1}+2}\subset
\cdots \subset S=S\cap F_{\alpha_1+d})\quad.$$
In particular, $\pi_\ualpha$ is a birational isomorphism. 

The varieties $\Vbb(\ualpha)$ are easily seen to be nonsingular,
as they may be realized as a tower of projective bundles over a point.

For a thorough treatment of Bott-Samelson varieties (in a more general
context) we refer the reader to \cite{math.AG/0302294}. We provide
here a self-contained construction of these varieties, adapted to
our application. In particular, we stress the independence of the 
definition of $\Vbb(\ualpha)$ from the ambient Grassmannian; we 
verify that the complement of $\Sbb(\ualpha)^\circ$ in $\Vbb(\ualpha)$ 
is a divisor with simple normal crossings; and we compute the 
push-forward ${\pi_{\ualpha}}_*$ explicitly.

\subsection{$\Vbb(\ualpha)$ as a tower of projective bundles}\label{secdesc}
The partition $\ualpha$ determines the following inductive
construction of $\Vbb(\ualpha)$. Each $\Vbb(\ualpha)$ will be a
nonsingular projective variety, of dimension $\sum_{i\ge 1} \alpha_i$,
endowed with bundles $\cL_i$, $\cQ_i$ of rank resp.~1,~$\alpha_i$. It
will map birationally to the corresponding Schubert variety
$\Sbb(\ualpha)$.

Here is the construction:
\begin{itemize}
\item For $\ualpha=(0)$, $\Vbb(\ualpha)$ is a point; for all $i\ge 0$, 
$\cQ_i=0$ and $\cL_i=\cO$;
\item Given any $\ualpha=(\alpha_1\ge \alpha_2\ge\dots)$, denote
by $\ualpha'$ the (shifted) `truncation' $(\alpha_2\ge \alpha_3\ge \dots)$,
and by $\cL'_i$, $\cQ'_i$ the corresponding bundles. Then
$\Vbb(\ualpha)$ is the following projective bundle over $\Vbb(\ualpha')$:
$$\xymatrix{
\Vbb(\ualpha):=\Pbb(\cQ'_1\oplus \cO^{\oplus(\alpha_1-\alpha_2+1)})
\ar[d]^{\rho_{\ualpha}^1} \\
\Vbb(\ualpha')
}$$
\item For $i\ge 2$, the distinguished bundles on $\Vbb(\ualpha)$ are
$$\cL_i:=(\rho_{\ualpha}^1)^*(\cL'_{i-1})\quad,\quad
\cQ_i:=(\rho_{\ualpha}^1)^*(\cQ'_{i-1})\quad.$$
The bundles $\cL_1$, $\cQ_1$ are defined by the tautological
sequence{\small
$$\xymatrix{
0 \ar[r] & \cL_1=\cO(-1) \ar[r] &
(\rho_{\ualpha}^1)^*(\cQ'_1\oplus \cO^{\oplus(\alpha_1-\alpha_2+1)})
= \cQ_2\oplus \cO^{\oplus(\alpha_1-\alpha_2+1)} \ar[r] &
\cQ_1 \ar[r] & 0
}$$}
\end{itemize}

\subsection{}

We let
$$\xi_i:=-c_1(\cL_i)\quad.$$
By construction, if $\alpha_j=0$ then $\cL_j\cong \cO$, 
and hence $\xi_j=0$. We record the following immediate consequence 
of the construction:

\begin{lemma}
Let $d$ be such that $\alpha_j=0$ for $j>d$. Then the Chow ring of 
$\Vbb(\ualpha)$ is generated by the classes $\xi_i$, $i=1,\dots,d$; thus, 
every element of $A_*\Vbb(\ualpha)$ may be expressed as 
$$\xi_1^{r_1}\cdots \xi_d^{r_d}\cap [\Vbb(\ualpha)]\quad,$$
for $r_1,\dots,r_d\ge 0$.
\end{lemma}

Relations among the monomials $\xi_1^{r_1}\cdots \xi_d^{r_d}$ follow
from the standard description of the Chow group of a projective bundle
(cf.~\cite{MR85k:14004}, \S3.3). The following observation will be crucial 
for our application, as it will be key to the computation of 
push-forwards (\S\ref{pf}), and is less immediate. We postpone its 
proof to \S\ref{proofofr}, at the end of this section.

We denote by $(N^d)$ the sequence
$$(\underbrace{N\ge \cdots\ge N}_d\ge 0 \ge\cdots)\quad,$$
corresponding to an $N\times d$ rectangle:
\begin{center}
\includegraphics[scale=.6]{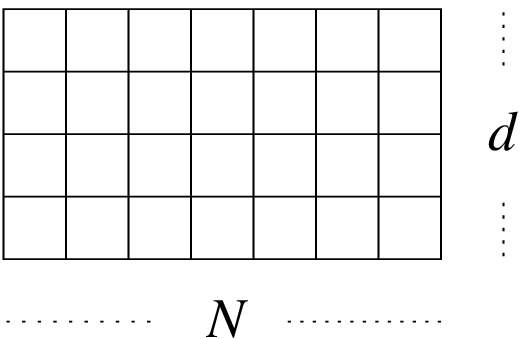}
\end{center}
(The corresponding Schubert variety $\Sbb((N^d))$ is then
isomorphic to $G_d(k^{N+d})$.)
Also, we denote by $\int$ the {\em degree\/} of a rational
equivalence class.

\begin{lemma}\label{reln}
Let $d\ge 2$, and let $r_1,\dots,r_d$ be nonnegative integers. 
Then for all $1\le i<d$ such that $r_i>0$:
$$\int (\prod_{j=1}^d \xi_j^{r_j})\cap [\Vbb((N^d))]=
-\int (\xi_1^{r_1}\cdots \xi_i^{r_{i+1}+1}
\xi_{i+1}^{r_i-1}\cdots \xi_d^{r_d})\cap[\Vbb((N^d))]\quad.$$
\end{lemma}

\subsection{}\label{strmaps}
By definition, if $\ubeta$ is obtained from $\ualpha$ by deleting the first
several entries, then there is a fibration
$$\xymatrix{
\rho:\Vbb(\ualpha) \ar@{->>}[r] & \Vbb(\ubeta)\quad,
}$$
which is in fact a composition of projective bundles; the
distinguished bundles on $\Vbb(\ubeta)$ pull back to distinguished
bundles with shifted indices on $\Vbb(\ualpha)$.
\begin{center}
\includegraphics[scale=.6]{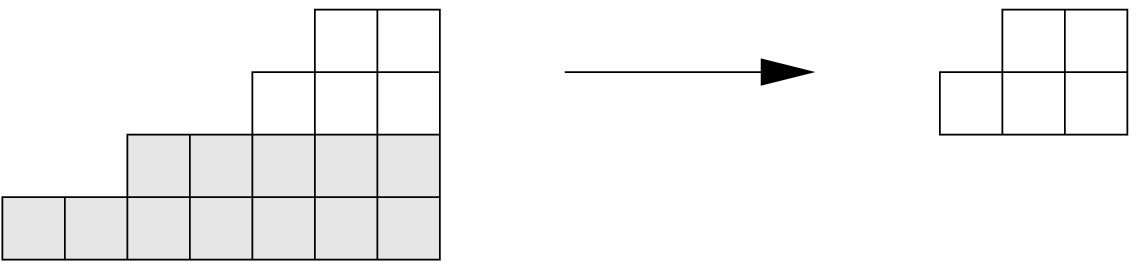}
\end{center}
If $\ubeta\le \ualpha$, then there are closed embeddings
$$\xymatrix{
\iota:\Vbb(\ubeta) \ar@{^(->}[r] & \Vbb(\ualpha)\quad,
}$$
such that (with evident notation)
$$\iota^*\cL_{\ualpha,i}=\cL_{\ubeta,i}\quad,\quad
\iota^*\cQ_{\ualpha,i}=\cQ_{\ubeta,i}\oplus\cO^{\alpha_i-\beta_i}\quad.$$
Indeed, assuming $\iota$ has been constructed for the truncated
sequences:
$$\xymatrix{
\iota:\Vbb(\ubeta') \ar@{^(->}[r] & \Vbb(\ualpha')\quad,
}$$
then the identification
$$\cQ'_{\ubeta,1}\oplus \cO^{\alpha_2-\beta_2}=\iota^*(\cQ'_{\ualpha,1})$$
gives an inclusion
$$\cQ'_{\ubeta,1}\oplus \cO^{\beta_1-\beta_2+1}
\subset 
\cQ'_{\ubeta,1}\oplus \cO^{\alpha_1-\beta_2+1}=
\cQ'_{\ubeta,1}\oplus \cO^{\alpha_2-\beta_2}\oplus \cO^{\alpha_1
-\alpha_2+1}=\iota^*(\cQ'_{\ualpha,1}\oplus \cO^{\alpha_1-\alpha_2+1})$$
projectivizing which induces the needed embedding $\Vbb(\ubeta) \subset
\Vbb(\ualpha)$.
\begin{center}
\includegraphics[scale=.6]{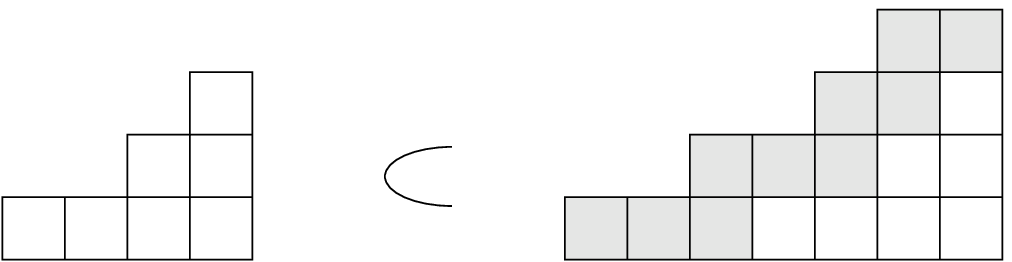}
\end{center}

\begin{lemma}\label{geomcla}
$[\Vbb(\ubeta)]=\xi_1^{\alpha_1-\beta_1}\cdots \xi_d^{\alpha_d-\beta_d}
\cap [\Vbb(\ualpha)]\in A_*(\Vbb(\ualpha))$.
\end{lemma}

\begin{proof}
This follows from a simple inductive argument, comparing the
sequences of relative tangent bundles for $\Vbb(\ualpha)$,
$\Vbb(\ubeta)$ over the varieties $\Vbb(\ualpha')$, $\Vbb(\ubeta')$ 
corresponding to the truncations (cf.~\cite{MR85k:14004}, \S B.5.8).
\end{proof}

Note that Lemma~\ref{geomcla} gives an explicit geometric realization
of monomials $\xi_1^{r_1}\cdots \xi_d^{r_d}$, for exponents $r_i$
such that the sequence $\alpha_1-r_1,\dots,\alpha_d-r_d$ is a
nonincreasing sequence of nonnegative integers.

\subsection{}\label{mapstograss}
The varieties $\Vbb(\ualpha)$ constructed above have natural maps
to Grassmannians. As in \S\ref{notat}, choose $d$ and $N$ so that
$$\left\{\aligned
&\text{$\alpha_i=0$ for $i>d$}\\
&\text{$N\ge \alpha_1$}
\endaligned
\right.\quad;$$
and fix a vector space $V$ of dimension~$N+d$, and a complete flag
$F_\bullet$ in $V$. These choices determine a map
$$\Vbb(\ualpha) \rightarrow G_d(V)$$
to the Grassmannian of $d$-planes in $V$, as follows. 

Denote by $\cF_r$ the trivial bundle with fiber $F_r$.

\begin{lemma}\label{epi}
For $i=1,\dots, d$ there are compatible epimorphisms
$$\cF_{\alpha_i+(d+1-i)} \twoheadrightarrow \cQ_i$$
of bundles over $\Vbb(\ualpha)$.
\end{lemma}

\begin{proof}
If $\ualpha$ is $0\ge 0\ge\dots$ there is nothing to show, for all $d$. 
For arbitrary $\ualpha$, assume that the epimorphisms have been 
constructed on the truncation $\Vbb(\ualpha')$; then the corresponding 
epimorphisms for $i\ge 2$ are obtained on $\Vbb(\ualpha)$ by pulling 
back via~$\rho_\ualpha^1$. In particular we have an epimorphism
$$\cF_{\alpha_2+(d-1)} \twoheadrightarrow \cQ_2\quad;$$
to obtain it for $i=1$:
$$\cF_{\alpha_1+d}=\cF_{\alpha_2+(d-1)}\oplus 
\cO^{\oplus(\alpha_1-\alpha_2+1)}\twoheadrightarrow
\cQ_2 \oplus \cO^{\oplus(\alpha_1-\alpha_2+1)}\twoheadrightarrow
\cQ_1\quad,$$
where the rightmost epimorphism comes from the universal 
sequence defining $\cQ_1$ on~$\Vbb(\ualpha)$.
\end{proof}

Let $\cS_i$ denote the kernel of the epimorphisms obtained in 
Lemma~\ref{epi}:
\begin{equation*}
\tag{*}
\xymatrix{
0 \ar[r] & \cS_i \ar[r] & \cF_{\alpha_i+(d+1-i)} \ar[r] & \cQ_i \ar[r]
& 0\quad;}
\end{equation*}
the rank of $\cS_i$ is $d+1-i$ (that is, the length of the $i$ row of the
`ladder'). By the universal property of Grassmannians, we obtain maps
$$\pi_\ualpha^i: \Vbb(\ualpha) \to G_{d+1-i}(F_{\alpha_i+(d+1-i)})$$
such that (*) is the pull-back of the 
tautological sequence over $G_{d+1-i}(F_{\alpha_i+(d+1-i)})$.

These maps are clearly compatible with the projections and
embeddings defined in~\S\ref{strmaps}. In fact,
a simple chase shows that the sequences (*) fit in diagrams 
with exact rows and columns:
$$\xymatrix{
& 0  & 0 & 0 & \\
0 \ar[r] & \cL_{i} \ar[r] \ar[u] & (\cQ_{i+1}\oplus 
\cO^{\oplus (\alpha_{i}-\alpha_{i+1}+1)}) \ar[r] \ar[u] & 
\cQ_{i} \ar[r] \ar[u] & 0\\
0 \ar[r] & \cS_{i} \ar[r] \ar@{->>}[u] & \cF_{\alpha_{i}+(d+1-i)} \ar[r] 
\ar@{->>}[u] & 
\cQ_{i} \ar[r] \ar@{=}[u] & 0 \\
0 \ar[r] & \cS_{i+1} \ar@{=}[r] \ar[u] & \cS_{i+1} \ar[r] \ar[u] & 
0 \ar[u] \\
& 0 \ar[u] & 0 \ar[u]
}$$
for all $i=1,\dots,d$.

\subsection{}\label{mainmapsec}
In particular, for $i=1$ we obtain a map
$$\pi_\ualpha: \Vbb(\ualpha) \to G_d(F_{\alpha_1+d})
\subset G_d(V)\quad,$$
where the last inclusion is induced by the inclusion
$F_{\alpha_1+d}\subset V$.

\begin{prop}\label{mainmap}
The map $\pi_\ualpha$ is a birational isomorphism onto~$\Sbb(\ualpha)$.
\end{prop}

Of course this is nothing but the map $\pi_\ualpha$ mentioned in
\S\ref{BSvars}: the fibers $\cS_i|_v$ over $v\in \Vbb(\ualpha)$ give 
a flag
$$\xymatrix{
\cS_d|_v \ar@{^(->}[r] \ar @{^(->}[d] & 
\cS_{d-1}|_v \ar@{^(->}[r] \ar @{^(->}[d] & 
\cdots  \ar@{^(->}[r] &
\cS_{1}|_v \ar @{^(->}[d] \\
F_{\alpha_d+1} \ar@{^(->}[r]  & 
F_{\alpha_{d-1}+2} \ar@{^(->}[r] & 
\cdots  \ar@{^(->}[r] &
F_{\alpha_1+d}
}$$
as prescribed in the description of $\Vbb(\ualpha)$ given in 
\S\ref{BSvars}, and $\pi_\ualpha(v)$ consists of the fiber
$\cS_1|_v$ viewed as a subspace of $V$. That is, $\pi_\ualpha(v)$
is a $d$-plane satisfying the incidence conditions defining 
$\Sbb(\ualpha)$. The image $\pi_\ualpha(v)$ is in the cell
$\Sbb(\ualpha)^\circ$ when
$$\dim (\cS_1|_v\cap F_{\alpha_i+(d+1-i)-1})<i\quad;$$
that is, when
$$\cS_i|_v\not\subset F_{\alpha_i+(d+1-i)-1}$$
for all $i=1,\dots,d$. The following proposition formalizes the
discussion given in \S\ref{BSvars}.

\begin{prop}\label{section}
The map $\pi_\ualpha$ admits a section over $\Sbb(\ualpha)^\circ$.
The complement of the image of $\Sbb(\ualpha)^\circ$ in $\Vbb(\ualpha)$ 
is a simple normal crossing divisor, whose components have class
$\xi_i$, $i\ge 1$.
\end{prop}

\begin{proof}
The section $i: \Sbb(\ualpha)^\circ \to \Vbb(\ualpha)$ may be defined
inductively. For $\ualpha=(0)$, both $\Sbb(\ualpha)^\circ$
and $\Vbb(\ualpha)$ are points. For arbitrary $\ualpha$, let
$\ualpha'=(\alpha_2\ge\alpha_3\ge\cdots)$ as usual, and assume
$$\xymatrix@M=10pt{
i'\,:\, \Sbb(\ualpha')^\circ \ar@{^(->}[r] & \Vbb(\ualpha')
}$$
has been constructed. For $S\in \Sbb(\ualpha)^\circ$, the intersection
of the corresponding $d$-plane with $F_{\alpha_{2}+(d-1)}$ has 
dimension exactly $d-1$, hence it determines a point $S'\in 
\Sbb(\ualpha')^\circ$.
$$\xymatrix@M=10pt{
\Sbb(\ualpha)^\circ \ar@{->>}[d] \ar@{^(..>}[r]^i & \Vbb(\ualpha) 
\ar@{->>}[d]^\rho \\
\Sbb(\ualpha')^\circ \ar@{^(->}[r]^{i'} & \Vbb(\ualpha')
}$$
By construction, $S'$ is naturally identified with the fiber 
$\cS_2|_{v'}$ of $\cS_2$ over $v'=i'(S')$. The one dimensional 
quotient $S/S'$ determines a one-dimensional subspace of
$$\frac{F_{\alpha_1+d}}{S'}=
\left(\frac{\cF_{\alpha_1+d}}{\cS_2}\right)|_{v'}=
(\cQ_2\oplus \cO^{\oplus \alpha_1-\alpha_2+1})|_{v'}\quad,$$
that is, a point $v$ of
$$\Pbb(\cQ_2\oplus \cO^{\oplus \alpha_1-\alpha_2+1})=\Vbb(\ualpha)$$
lying over $v'$. Setting $i(S)=v$ lifts $i'$, as needed.

The statement on the complement of $\Sbb(\ualpha)^\circ$ in
$\Vbb(\ualpha)$ may also be verified inductively. The complement
$D'$ of $\Sbb(\ualpha')^\circ$ in $\Vbb(\ualpha')$ may be assumed
to be a simple normal crossing divisor, with components of class
$\xi_i$, $i\ge 2$. For $S$ over $S'\in \Sbb(\ualpha')^\circ$, the 
condition $S\in \Sbb(\ualpha)^\circ$ is equivalent to
$S\not\subset F_{\alpha_1+d-1}$,
that is (with notation as above) to
$$v\not\in \Pbb(\cQ_2\oplus \cO^{\oplus \alpha_1-\alpha_2})\quad.$$
Therefore, the complement of $\Sbb(\ualpha)^\circ$ in $\Vbb(\ualpha)$
consists of $D=E\cup \rho^{-1}(D')$, where $E:=\Pbb(\cQ_2\oplus 
\cO^{\oplus \alpha_1-\alpha_2})$ is a hypersurface of $\Vbb(\ualpha)$ 
of class $\xi_1$. It remains to verify that $D$ is a divisor with
simple normal crossings, and this is an instance of the following
more general situation: let $\pi:Y \to X$ be a smooth morphism of 
varieties, $E \subset Y$ an irreducible divisor, smooth over $X$, 
and $D'$ a simple normal crossing divisor on $X$. Then 
$\pi^{-1}(D')\cup E$ is a simple normal crossing divisor on $Y$,
as needed.
\end{proof}

\subsection{Pull-back}\label{pb}
The left column in the diagram displayed at the end of
\S\ref{mapstograss} is the exact sequence
$$\xymatrix{
0 \ar[r] & \cS_{i+1} \ar[r] & \cS_i \ar[r] & \cL_i \ar[r] & 0
}$$
($i=1,\dots,d$). The existence of these sequences amounts to
the fact that the variety $\Vbb(\ualpha)$ is a concrete
realization of the splitting construction for the restriction
of the tautological bundles to the Schubert varieties
$\Sbb(\alpha_i\ge \alpha_{i+1}\ge \dots)$: for each $i=1,\dots,d$
we have constructed a map
$$\xymatrix{
\pi_{\ualpha}^i\,:\, \Vbb(\ualpha) \ar@{->>}[r] & 
\Sbb(\alpha_i\ge \alpha_{i+1}\ge \dots) \subset G_{d+1-i}(V)
}$$
such that $\cS_i$ is the pull-back of the tautological subbundle.

\begin{lemma}\label{classS}
With notation as above:
\begin{itemize}
\item $c(\cS_i)=(1+\xi_i)\cdots (1+\xi_d)$;
\item the image of the pull-back map 
$${\pi_{\ualpha}^i}^*: A^*(G_{d+1-i}(V)) \to A^*(\Vbb(\ualpha))$$
consists of the symmetric functions in $\xi_i,\dots, \xi_d$.
\end{itemize}\end{lemma}

\begin{proof}
The first statement is an immediate consequence of the exact
sequences recalled above. The second follows from the fact that
the Chow ring of the Grassmannian is generated by the Chern
classes of the tautological subbundle, and the considerations
preceding the statement.
\end{proof}

\subsection{Push-forward}\label{pf}
Lemmas~\ref{reln} and~\ref{classS} allow us to determine the
push-forward
$${\pi_{\ualpha}}_*: A_*(\Vbb(\ualpha)) \to A_*(G_d(V))
\quad,$$
and this will be necessary for our application.

We need another piece of notation. Let $a_1, \dots,a_d$ denote
integers; if $a_1>\dots> a_d>0$, we let
$$\Omega(a_d,\dots,a_1)=\Sbb(a_1-d\ge a_2-(d-1)\ge\dots \ge 
a_d-1\ge 0\ge\dots)\quad.$$
The notation is chosen in order to match standard terminology
(cf.~\cite{MR85k:14004}, \S14.7).
The number $a_i$ denotes the number of boxes in the
$i$-th row (from the bottom) of the usual diagram, {\em with\/}
a $d$-ladder adjoined to the right:
\begin{center}
\includegraphics[scale=.6]{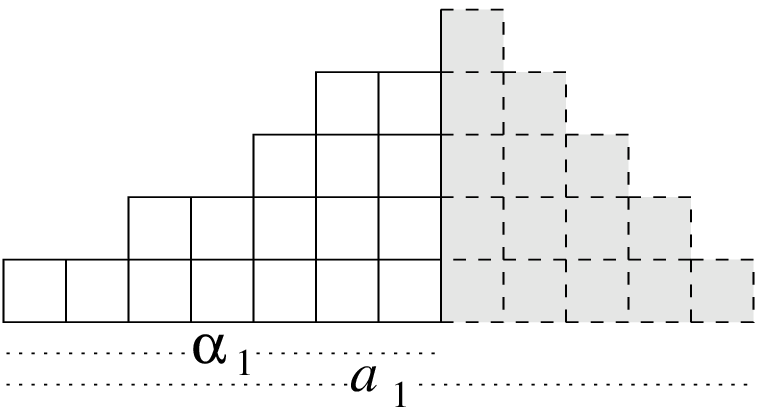}
\end{center}
For arbitrary $a_1,\dots,a_d$, define the class 
$[\Omega(a_d,\dots,a_1)]$ in the Chow group of a Schubert variety
to be
$$[\Omega(a_d,\dots,a_1)]=(-1)^\sigma\, [\Omega(a_{\sigma(d)},\dots,
a_{\sigma(1)})]$$
if $a_1,\dots,a_d$ are positive and distinct, and $\sigma\in S_d$ 
is the permutation such that $a_{\sigma(d)}< \cdots < a_{\sigma(1)}$;
and 0 if the integers $a_1,\dots,a_d$ are not positive and distinct.

\begin{prop}\label{pfprop}
Let $r_1,\dots,r_d$ be nonnegative integers. Then
$${\pi_{\ualpha}}_*\left(\xi_1^{r_1}\cdots \xi_d^{r_d}\cap
[\Vbb(\ualpha)]\right)=[\Omega(\alpha_d-r_d+1,\dots,\alpha_1-r_1+d)]
\quad.$$
\end{prop}

\begin{proof}
The formula is immediate if $d=1$, so we assume $d\ge 2$. 
By compatibility with embeddings (\S\ref{strmaps}), we can assume
that $\ualpha=(N^d)$, that is, $\alpha_1=\dots=\alpha_d=N$,
$\alpha_{d+1}=0$. 
If $0\le r_1\le \cdots \le r_d\le N$ then
$$\xi_1^{r_1}\cdots \xi_d^{r_d}\cap [\Vbb((N^d))]=
[\Vbb(\ubeta)]$$
for $\ubeta=(N-r_1\ge \cdots N-r_d\ge 0\cdots)$,
\begin{center}
\includegraphics[scale=.6]{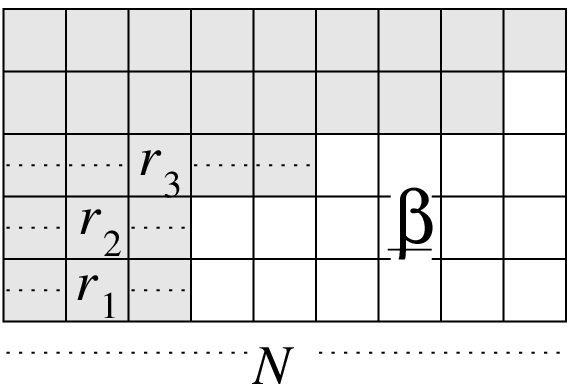}
\end{center}
by Lemma~\ref{geomcla}. Thus
$${\pi_\ualpha}_*\left(\xi_1^{r_1}\cdots \xi_d^{r_d}\cap 
[\Vbb((N^d))]\right)={\pi_\ualpha}_*\left([\Vbb(\ubeta)]
\right)=[\Sbb(\ubeta)]=[\Omega(N-r_d+1,\dots,N-r_1+d)]$$
if $0\le r_1\le \cdots \le r_d\le N$, that is, if
$$0<N-r_d+1<\cdots <N-r_1+d\quad.$$
In order to prove the formula in the general case, it suffices
to show that both sides behave in the same way after permutations
of this list of integers; that is, it suffices to show that
$${\pi_\ualpha}_*\left(\xi_1^{s_1}\cdots \xi_d^{s_d}\cap 
[\Vbb((N^d))]\right)=(-1)^\sigma 
{\pi_\ualpha}_*\left(\xi_1^{r_1}\cdots \xi_d^{r_d}\cap 
[\Vbb((N^d))]\right)$$
if the list
$$N-s_d+1\quad,\quad\dots\quad, \quad N-s_1+d\quad$$
is obtained by applying the permutation $\sigma$ to the list
$$N-r_d+1\quad,\quad\dots\quad, \quad N-r_1+d\quad.$$
As transpositions generate the group of permutations, it suffices
to prove that
$${\pi_\ualpha}_*\left(\xi_1^{r_1}\cdots \xi_i^{r_{i+1}+1}
\xi_{i+1}^{r_i-1}\cdots \xi_d^{r_d}\cap [\Vbb((N^d))]\right)
=-{\pi_\ualpha}_*\left(\xi_1^{r_1}\cdots \xi_d^{r_d}\cap 
[\Vbb((N^d))]\right)$$
for all nonnegative integer $r_1,\dots,r_d$, assuming $r_i>0$.
By Poincar\'e duality in the Grassmannian and the projection
formula, it suffices to prove that
$$\int\xi_1^{r_1}\cdots \xi_i^{r_{i+1}+1}\xi_{i+1}^{r_i-1}\cdots 
\xi_d^{r_d}\cdot \pi_\ualpha^*(C)\cap [\Vbb((N^d))]
=-\int\xi_1^{r_1}\cdots \xi_d^{r_d}\cdot \pi_\ualpha^*(C)
\cap [\Vbb((N^d))]$$
for all classes $C$ of codimension $dN-\sum r_j$
in the Grassmannian.  

By Lemma~\ref{classS}, 
$$\pi_\ualpha^*(C)=P(\xi_1,\dots,\xi_d)$$
for a homogeneous symmetric polynomial $P(\xi_1,\dots,\xi_d)$ of 
degree $dN-\sum r_j$; we have to prove that 
$$\xi_1^{r_1}\cdots \xi_d^{r_d}\cdot P(\xi_1,\dots,\xi_d)$$
and
$$-\xi_1^{r_1}\cdots \xi_i^{r_{i+1}+1}\xi_{i+1}^{r_i-1}\cdots 
\xi_d^{r_d}\cdot P(\xi_1,\dots,\xi_d)$$
agree in $[\Vbb((N^d))]$.
Now, there is a one-to-one correspondence between monomials
in these expressions: each monomial 
$$(\prod_{j=1}^d \xi_j^{r_j}) \cdots \xi_i^a \xi_{i+1}^b\cdots$$
in the first one corresponds to exactly one monomial
$$-(\xi_1^{r_1}\cdots \xi_i^{r_{i+1}+1}
\xi_{i+1}^{r_i-1}\cdots \xi_d^{r_d}) \cdots \xi_i^b \xi_{i+1}^a\cdots$$
in the second one; it suffices to show that these two terms match. 
After absorbing all exponents into the $r_i$'s, the stated equality
is then reduced to the statement of Lemma~\ref{reln} (proved in 
\S\ref{proofofr} below), and we are done.
\end{proof}

The push-forward operation may be visualized as follows:
adjoin the $d$-ladder to~$\ualpha$; take away boxes from the
rows as dictated by the given monomial in $\xi_1,\dots,\xi_d$; 
rearrange the remaining rows in {\em strictly} increasing order, 
keeping track of the sign of the needed permutation; the push-forward 
is then read from the complement of the ladder.

\begin{example}
Let $\ualpha=(7\ge 6\ge 4 \ge 3)$; then
$${\pi_{\ualpha}}_*\left(\xi_1^{10}\xi_2^2 \cap [\Vbb(\ualpha)]
\right)=
-[\Sbb(3\ge 3\ge 2)]\quad:$$
\begin{center}
\includegraphics[scale=.6]{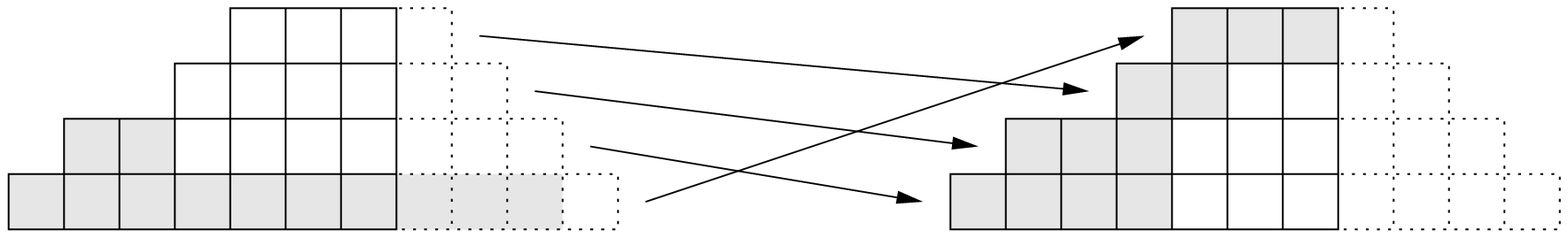}
\end{center}

The push-forward is~0 if these operations cannot be performed;
for example, if the numbers of `white' boxes in the rows are
not distinct:
$${\pi_{\ualpha}}_*\left(\xi_1^{10}\xi_2^3 \cap [\Vbb(\ualpha)]
\right)= 0\quad.$$
\begin{center}
\includegraphics[scale=.6]{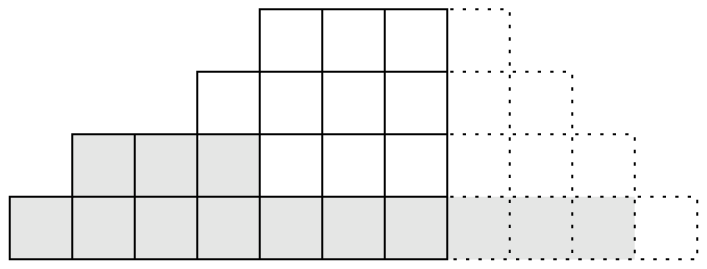}
\end{center}
\end{example}

\subsection{Proof of Lemma~\ref{reln}}\label{proofofr}
Let 
$$\rho: \Vbb((N^d)) \to \Vbb(0) = \text{{point}}$$
be the projection; Lemma~\ref{reln} is equivalent to
$$\rho_* \left(\xi_1^{r_1}\cdots \xi_d^{r_d}\cap [\Vbb((N^d))]
\right)=
-\rho_* \left(\xi_1^{r_1}\cdots \xi_i^{r_{i+1}+1}
\xi_{i+1}^{r_i-1}\cdots \xi_d^{r_d}\cap[\Vbb((N^d))]
\right)\quad.$$
This reduces the lemma to a particular case of the following 
template situation. 

Let $X$ be any scheme, and let $X_2=\Pbb(\cE)$ be a projective 
bundle over $X$, with $\rk\cE=N+1$. The universal sequence 
over $X_2$:
$$\xymatrix{
0 \ar[r] & \cO(-1) \ar[r] & \cE \ar[r] & \cQ \ar[r] & 0
}$$
determines a quotient bundle $\cQ$, of rank $N$. Let
$X_1=\Pbb(\cQ\oplus \cO)$, projecting to $X$:
$$\xymatrix{
\rho: X_1 \ar[r]^{\rho_1} & X_2\ar[r]^{\rho_2} \ar[r] & X
}\quad.$$
Lemma~\ref{reln} is an immediate consequence of the following 
explicit computation.

\begin{lemma}
Let $\xi_2=-c_1(\cO_{\Pbb(\cE)}(-1))$ on $X_2$ (and its pull-back 
to $X_1$); and let $\xi_1=-c_1(\cO_{\Pbb(\cQ\oplus \cO)})(-1))$ on 
$X_1$. Then
$$\rho_*(\xi_1^a \xi_2^b\cap [X_1])
=-\rho_*(\xi_1^{b+1} \xi_2^{a-1}\cap [X_1])\quad.$$
\end{lemma}

\begin{proof}
By the very definition of Segre class,
$${\rho_1}_*(\sum_j \xi_1^j\cap [X_1])=s(\cQ\oplus \cO)\cap [X_2]
=s(\cE)c(\cO_{\Pbb(\cE)}(-1))\cap [X_2]=s(\cE)(1-\xi_2)\cap [X_2]\quad.$$
Therefore
$${\rho_1}_*(\xi_1^a\cap [X_1])=(s_{a-N}(\cE)-s_{a-1-N}(\cE)
\,\xi_2)\cap [X_2]\quad.$$
By the same token,
$${\rho_2}_*(\sum_j \xi_2^j \cap [X_2])=s(\cE)\cap [X]\quad,$$
and it follows that
$$\rho_*(\xi_1^a\xi_2^b\cap [X_1])
=(s_{a-N}(\cE)s_{b-N}(\cE)-s_{a-1-N}(\cE)s_{b+1-N}(\cE))
\cap [X]\quad.$$
The substitution $a\mapsto b+1$, $b\mapsto a-1$ switches the summands
in the right-hand-side, proving the statement.
\end{proof}


\section{Chern-Schwartz-MacPherson classes of Schubert cells}\label{CSM}

\subsection{Chern classes of singular varieties}
We now assume that the ground field $k$ has characteristic~0.  As
recalled in the introduction, there is a good theory of Chern classes
for singular varieties: $\csm(X)\in A_*(X)$ will denote the {\em
Chern-Schwartz-MacPherson\/} class of the (possibly) singular variety
$X$. There are a number of different approaches to the definition of
these classes; the following is best suited to our purposes. A more
thorough discussion may be found in \cite{math.AG/0507029} or
\cite{MR2209219}.

Let $X$ be a complete variety over an algebraically closed field of
characteristic zero. We denote by $F(X)$ the group of {\em
constructible functions\/} on $X$, that is, the free abelian group on
characteristic functions of subvarieties of $X$. Thus, every
constructible function may be written (uniquely) as $\sum m_Z
\,\one_Z$, where the sum ranges over finitely many closed subvarieties
$Z\subset X$, $m_Z\in \Zbb$, and $\one_Z$ denotes the function with
value $1$ on $Z$ and $0$ on its complement.

Constructible functions may also be written as finite sums $\sum m_W
\one_W$, for $W\subset X$ {\em locally\/} closed subvarieties of $X$;
of course these expressions are no longer unique. For a given
$\varphi\in F(X)$, we consider any decomposition
$$X=\bigcup_{j\in J} W_j$$
of $X$ into finitely many, locally closed, {\em nonsingular\/}
subvarieties $W_j$, such that
$$\varphi=\sum_{j\in J} m_j\, \one_{W_j}\quad,$$
for $m_j\in \Zbb$. We use such a decomposition to associate to
$\varphi$ a class
$$c_*(\varphi)\in A_*X\quad,$$
as follows.

By resolution of singularities, if $W\subset X$ is nonsingular then
there exists a nonsingular completion $\overline W$ of $W$, and 
a proper morphism 
$$\omega: \overline W \to X$$
such that the complement $\overline W\setminus W$ is a divisor $D$
with normal crossings and nonsingular components $D_i$, $i\in I$.

\begin{defin}\label{defcsm}
For $W\subset X$ a nonsingular locally closed subvariety, we set
$$c_*(\one_W):=\omega_*\left(c(T\overline W(-\log D))
\cap [\overline W]\right)\in A_*(X)\quad.$$
\end{defin}

Here $T\overline W(-\log D)$ denotes the bundle of vector fields
with logarithmic zeros along the components of $D$; as is well
known,
$$c(T\overline W(-\log D))=\frac{c(T\overline W)}{\prod_{i\in I}
(1+D_i)}\quad.$$
One can verify that the class $c_*(\one_W)$ is independent of the 
chosen completion~$\overline W$.

The class $c_*(\varphi)$ is defined by linearity: for a decomposition
$X=\cup_{j\in J} W_j$ as above, we set
$$c_*(\varphi):=\sum_{j\in J} m_j\, c_*(\one_{W_j})\quad.$$

One can verify that the class $c_*(\varphi)$ only depends on the 
constructible function $\varphi$, not on the chosen decomposition. 
In fact (cf.~\cite{MR2209219}, Th\'eor\`eme~3.3), the homomorphism
$c_*: F(X) \to A_*(X)$ agrees with the one induced by MacPherson's 
natural transformation
from the functor~$F$ (with covariance defined by Euler characteristic
of fibers, cf.~\cite{MR50:13587} or \cite{MR85k:14004}, Example~19.1.7) 
to the Chow group functor $A_*$. In particular, for $\varphi=\one_X$ 
one obtains the {\em Chern-Schwartz-MacPherson\/}
class of~$X$:
$$\csm(X):= c_*(\one_X)\in A_*(X)\quad;$$
and the preceding discussion shows that if
$$X=\coprod_{j\in J} W_j$$
is a decomposition of $X$ into disjoint nonsingular locally closed 
subvarieties, then
$$\csm(X)=\sum_{j\in J}\,c_*(\one_{W_j})\quad.$$
We denote $c_*(\one_{W_j})$ by $\csm(W_j)$, for notational consistency;
this abuse of language is harmless in context, provided that the reader 
keeps in mind that $\csm(W_j)$ is a class in the Chow group of the 
ambient variety $X$, not of $W_j$.

\subsection{}
We are essentially ready to compute $\csm(\Sbb(\ualpha)^\circ)$, by 
applying the formula given in Definition~\ref{defcsm}; the only missing 
ingredient is the Chern class of the tangent bundle of $\Vbb(\ualpha)$. 

Recall that $\cQ_i$ is trivial for $i\gg 0$; thus 
$c(\cL_i^\vee \otimes\cQ_i)=1$ for $i\gg 0$, and only finitely many terms 
contribute nontrivially to the formula in the following statement.

\begin{prop}\label{tang}
$$c(T\Vbb(\ualpha))=\prod_{i\ge 1} c(\cL_i^\vee\otimes\cQ_i)
\quad.$$
\end{prop}

\begin{proof}
As $\Vbb(\ualpha)$ is defined inductively, we have to check that this
is correct for $\ualpha=(0)$, which is trivially the case, and that
the classes have identical behavior when going from the truncation
$\ualpha'=(\alpha_2\ge \dots)$ to $\ualpha=(\alpha_1\ge \alpha_2\ge
\dots)$.

Now, by definition,
$$\Vbb(\ualpha)=\Pbb(\cQ_2\oplus \cO^{\oplus \alpha_1-\alpha_2+1})$$
over $\Vbb(\ualpha')$, 
with relative tangent bundle $\cL_1^\vee\otimes\cQ_1$: this is computed
by tensoring the universal sequence
$$\xymatrix{
0 \ar[r] & \cL_1 \ar[r] & \cQ_2\oplus \cO^{\oplus \alpha_1-\alpha_2+1}
\ar[r] & \cQ_1 \ar[r] & 0
}$$
by $\cL_1^\vee$ (cf.~\cite{MR85k:14004}, \S B.5.8.\/)
Then the statement follows immediately.
\end{proof}

The main results of this paper will be obtained from the following
consequence.

\begin{corol}\label{urform}
Let $\pi_\ualpha:\Vbb(\ualpha) \to \Sbb(\ualpha)$ be the birational 
isomorphism defined in \S\ref{mainmapsec}.

Then the Chern-Schwartz-MacPherson class of the Schubert cell 
$\Sbb(\ualpha)^\circ$ is
\begin{align*}
\csm(\Sbb(\ualpha)^\circ)&=
{\pi_{\ualpha}}_*\left(
\frac{c(T\Vbb(\ualpha))}{\prod_{i\ge 1} c(\cL_i^\vee)}
\cap[\Vbb(\ualpha)]\right)\\
&={\pi_{\ualpha}}_*\left(\prod_{i\ge 1}
c(\cL_i^\vee)^{\alpha_i-\alpha_{i+1}}\, c(\cL_i^\vee \otimes \cQ_{i+1})
\cap[\Vbb(\ualpha)]\right)
\quad.
\end{align*}
\end{corol}

\begin{proof}
By Proposition~\ref{section}, we can identify $\Sbb(\ualpha)^\circ$ with 
an open dense subset of $\Vbb(\ualpha)$, and the
complement of $\Sbb(\ualpha)^\circ$ in $\Vbb(\ualpha)$ is a divisor with 
normal crossings and nonsingular components, with class 
$\xi_i=c_1(\cL_i^\vee)$.  According to Definition~\ref{defcsm},
$$\csm(\Sbb(\ualpha)^\circ)={\pi_{\ualpha}}_*\left(
\frac{c(T\Vbb(\ualpha))}{\prod_{i\ge 1} c(\cL_i^\vee)}\cap[\Vbb(\ualpha)]
\right)$$
(since $\cL_i=\cO$ for $i\gg 0$, only finitely many terms contribute
to the denominator).

Now tensor the sequences
$$\xymatrix{
0 \ar[r] & \cL_i \ar[r] & \cQ_{i+1}\oplus \cO^{\oplus \alpha_i-\alpha_{i+1}+1}
\ar[r] & \cQ_i \ar[r] & 0
}$$
by $\cL_i^\vee$, obtaining
$$c(\cL_i^\vee\otimes\cQ_i)=c(\cL_i^\vee\otimes \cQ_{i+1})\,
c(\cL_i^\vee)^{\alpha_i-\alpha_{i+1}+1}\quad;$$
applying Proposition~\ref{tang} gives
$$\csm(\Sbb(\ualpha)^\circ)={\pi_{\ualpha}}_*\left(\prod_{i\ge 1}
\frac{c(\cL_i^\vee\otimes \cQ_{i+1})\, c(\cL_i^\vee)^{\alpha_i-\alpha_{i+1}+1}}
{c(\cL_i^\vee)}\cap[\Vbb(\ualpha)]\right)$$
with the stated result.
\end{proof}

The next task is to obtain explicit formulas from the statement of 
Corollary~\ref{urform}. We offer several versions, all easily amenable to 
computer implementation. Since $\csm(\Sbb(\ualpha)^\circ)\in 
A_*\Sbb(\ualpha)$, 
$$\csm(\Sbb(\ualpha)^\circ) = \sum_{\ubeta\le \ualpha} \gamma_{\ualpha,
\ubeta} [\Sbb(\ubeta)]\quad,$$
for uniquely determined coefficients $\gamma_{\ualpha,\ubeta}\in \Zbb$; 
we give formulas computing these coefficients.

The reader should note that the formulas will often appear to depend 
on the choice of an integer $d$ such that $\alpha_{d+1}=0$; 
the results of the computations must be independent of this choice,
since so is the class $\csm(\Sbb(\ualpha)^\circ)$.

\subsection{Explicit computations of $\gamma_{\ualpha,\ubeta}$}
We denote by 
$$h_a(y_1,\dots,y_n)=\sum_{i_1+\cdots +i_n=a} y_1^{i_1}\cdots y_n^{i_n}$$
the {\em complete symmetric polynomial\/} of degree $a$ in the variables 
$y_1,\dots,y_n$.

\begin{theorem}\label{complsym}
Let $\ualpha=(\alpha_1\ge\dots\ge \alpha_d)$, $\ubeta=(\beta_1
\ge\dots\ge\beta_d)$ be partitions.

For $b_1,\dots,b_d$ positive integers, let 
$C_{\ualpha}(b_1,\dots,b_d)$ be the coefficient of 
$$x_1^{\alpha_1+d-b_1}\cdots x_d^{\alpha_d+1-b_d}$$
in the polynomial
$$\prod_{i=1}^d (1+x_i)^{\alpha_i-\alpha_{i+1}} \cdot
h_{\alpha_{i+1}}(1+x_i,x_{i+1},\dots, x_d)\quad.$$

Then
$$\gamma_{\ualpha,\ubeta}= \sum_{\sigma\in S_d} (-1)^\sigma
C_\ualpha (b_{\sigma(1)},\dots,b_{\sigma(d)})\quad,$$
where $b_i=\beta_i+(d+1-i)$.
\end{theorem}

\begin{proof}
By Corollary~\ref{urform},
$$\csm(\Sbb(\ualpha)^\circ)={\pi_{\ualpha}}_*\left(\prod_{i=1}^d
c(\cL_i^\vee)^{\alpha_i-\alpha_{i+1}}\, c(\cL_i^\vee \otimes \cQ_{i+1})
\cap[\Vbb(\ualpha)]\right)\quad,$$
where $d$ is such that $\alpha_{d+1}=0$ (so that $\cL_i=\cO$
and $\cQ_i=0$ for $i>d$).

Recall that for all $i$ there are exact sequences
((*) in \S\ref{mapstograss})
$$\xymatrix{
0 \ar[r] & \cS_i \ar[r] & \cF_{\alpha_i+(d+1-i)} \ar[r] & \cQ_i \ar[r]
& 0\quad;
}$$
and in particular
$$c(\cQ_{i+1})=\frac 1{c(\cS_{i+1})} = 
\frac 1{(1-\xi_{i+1})\cdots (1-\xi_d)}
=\sum_{j=1}^{\alpha_{i+1}} h_j(\xi_{i+1},\dots,\xi_d)\quad.$$
Therefore (\cite{MR85k:14004}, Example~3.2.2)
$$c(\cL_i^\vee\otimes \cQ_{i+1})=h_{\alpha_{i+1}} (1+\xi_i,\xi_{i+1},
\dots, \xi_d)$$
and hence
$$\prod_{i=1}^d c(\cL_i^\vee)^{\alpha_i-\alpha_{i+1}}\, 
c(\cL_i^\vee \otimes \cQ_{i+1})=
\prod_{i=1}^d (1+\xi_i)^{\alpha_i-\alpha_{i+1}}\, 
h_{\alpha_{i+1}} (1+\xi_i,\xi_{i+1}, \dots, \xi_d)\quad.$$
Now let $x_1,\dots,x_d$ be variables, and write
$$\prod_{i=1}^d (1+x_i)^{\alpha_i-\alpha_{i+1}}\, 
h_{\alpha_{i+1}} (1+x_i,x_{i+1}, \dots, x_d)
=\sum_{r_1,\dots,r_d\ge 0} e_{r_1,\dots,r_d}\, x_1^{r_1} \cdots
x_d^{r_d}\quad;$$
by Proposition~\ref{pfprop},
\begin{align*}
\csm(\Sbb(\ualpha)^\circ) &= {\pi_{\ualpha}}_*\left(
\sum_{r_1,\dots,r_d\ge 0} e_{r_1,\dots,r_d}\, \xi_1^{r_1} \cdots
\xi_d^{r_d}\cap [\Vbb(\ualpha)]\right)\\
& =\sum_{r_1,\dots,r_d\ge 0} e_{r_1,\dots,r_d}\, 
[\Omega(\alpha_d-r_d+1,\dots,\alpha_1-r_1+d)]
\quad,
\end{align*}
with notation as in \S\ref{pf}. By definition, 
$\gamma_{\ualpha,\ubeta}$ is the coefficient of $[\Sbb(\beta)]$ in
this expression.
Let then $b_1=\beta_1+d$, \dots, $b_d=\beta_d+1$; thus $b_1>\cdots
>b_d$ are positive integers, $[\Sbb(\ubeta)]=[\Omega(b_d,\dots,b_1)]$, 
and
$$[\Omega(b_{\sigma(d)},\dots,b_{\sigma(1)})]=(-1)^\sigma
[\Sbb(\ubeta)]\quad.$$
With this notation, the coefficient of $[\Sbb(\ubeta)]$ is given by
$$\sum_{\sigma\in S_d} (-1)^\sigma e_{\alpha_1+d-b_{\sigma(1)},
\dots, \alpha_d+1-b_{\sigma(d)}}=
\sum_{\sigma\in S_d} (-1)^\sigma C_\ualpha(b_{\sigma(1)},
\dots, b_{\sigma(d)})\quad,$$
as stated.
\end{proof}

\begin{example}\label{exa22}
For two-row diagrams, i.e., $\ualpha=(\alpha_1\ge \alpha_2)$, the 
CSM class of the corresponding Schubert cell may be 
read off the expansion of
$$(1+x_1)^{\alpha_1-\alpha_2}(1+x_2)^{\alpha_2}
\sum_{i=0}^{\alpha_2}(1+x_1)^i x_2^{\alpha_2-i}\quad.$$
For the open cell in $G_2(\Cbb^4)$ (=the Grassmannian of 
projective lines in $\Pbb^3$):
\begin{center}
\includegraphics[scale=.6]{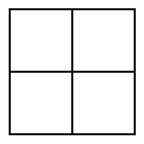}
\end{center}
this is
$$1+3\,x_2+ 4\,x_2^2+ 5\,x_1x_2+ x_1^2+4\,x_1 x_2^2+x_1^2 x_2^2
+\left(2\,x_1+2\,x_1^2 x_2+3\,x_2^3+ x_1x_2^3+ x_2^4\right)\quad.$$
The coefficient of the Schubert class corresponding to
\begin{center}
\includegraphics[scale=.6]{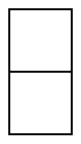}
\end{center}
that is, of $[\Sbb(1\ge 1)]$, is
$$\gamma_{(2\ge 2),(1\ge 1)}=C_{2\ge 2}(3,2)-C_{2\ge 2}(2,3)=
e_{1,1}-e_{2,0}=5-1=4\quad.$$
The terms in parentheses do not correspond to any subdiagram of
$\ualpha$; they do correspond to classes in $\Vbb(\ualpha)$ 
(setting $x_i=\xi_i$), but they all push-forward to $0$ in 
$A_*(G_2(\Cbb^4))$, as may be checked by applying 
Proposition~\ref{pfprop}.
\end{example}

There is a useful alternative to the expression given in 
Theorem~\ref{complsym}:

\begin{theorem}\label{ratf}
Let $\ualpha=(\alpha_1\ge\dots\ge \alpha_d)$, $\ubeta=(\beta_1
\ge\dots\ge\beta_d)$ be partitions.

For $b_1,\dots,b_d$ positive integers, let 
$C'_{\ualpha}(b_1,\dots,b_d)$ be the coefficient of 
$$x_1^{\alpha_1+d-b_1}\cdots x_d^{\alpha_d+1-b_d}$$
in the expansion of the rational function
$$\frac{\prod_{i=1}^d (1+x_i)^{\alpha_i+(d-i)}}
{\prod_{1\le i< j\le d} (1+x_i-x_j)}$$
at $0$. Then 
$$\gamma_{\ualpha,\ubeta}= \sum_{\sigma\in S_d} (-1)^\sigma
C'_\ualpha (b_{\sigma(1)},\dots,b_{\sigma(d)})\quad,$$
where $b_i=\beta_i+(d+1-i)$.
\end{theorem}

\begin{proof}
Use again the exact sequence
$$\xymatrix{
0 \ar[r] & \cS_i \ar[r] & \cF_{\alpha_i+(d+1-i)} \ar[r] & \cQ_i \ar[r]
& 0
}$$
to obtain
$$c(\cL_i^\vee\otimes \cQ_i)=\frac{c(\cL_i^\vee)^{\alpha_i+(d+1-i)}}
{c(\cL_i^\vee \otimes\cS_i)}=\frac{c(\cL_i^\vee)^{\alpha_i+(d+1-i)}}
{\prod_{i\le j\le d} c(\cL_i^\vee \otimes\cL_j)}\quad,$$
and hence
$$\frac{c(T\Vbb(\ualpha))}{\prod_{1\le i\le d} c(\cL_i^\vee)}
=\prod_{1\le i\le d} \frac{c(\cL_i^\vee)^{\alpha_i+(d-i)}}
{\prod_{i< j\le d} c(\cL_i^\vee \otimes\cL_j)}
=\frac{\prod_{1\le i\le d}(1+\xi_i)^{\alpha_i+(d-i)}}
{\prod_{1\le i<j\le d}(1+\xi_i-\xi_j)}\quad.$$
Then argue as in the proof of Theorem~\ref{complsym}.
\end{proof}

\begin{example}
The numbers $C_{\ualpha}(b_1,\dots,b_d)$ and $C'_{\ualpha}(b_1,\dots,b_d)$ 
are {\em not\/} equal in general. For example,
$$C_{(2\ge 1\ge 1)}(5,2,1)=4 \quad,\quad C'_{(2\ge 1\ge 1)}(5,2,1)=5
\quad.$$
It is a consequence of Theorem~\ref{complsym} and~\ref{ratf} that,
however, the `antisymmetrization' of these coefficients must agree.
Thus
$$C_{(2\ge 1\ge 1)}(5,2,1)-C_{(2\ge 1\ge 1)}(5,1,2) = 3
=C'_{(2\ge 1\ge 1)}(5,2,1)-C'_{(2\ge 1\ge 1)}(5,1,2)$$
both compute the coefficient $\gamma_{(2\ge 1\ge 1),(2)}$ of the 
Schubert class corresponding to
\begin{center}
\includegraphics[scale=.6]{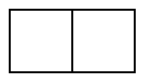}
\end{center}
in the Chern-Schwartz-MacPherson class of the Schubert cell of
\begin{center}
\includegraphics[scale=.6]{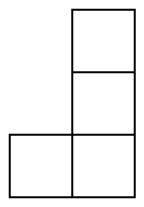}
\end{center}
(the other permutations of the arguments $5,2,1$ all give 
vanishing contributions to the computation of $\gamma_{(2\ge 1\ge 1),
(2\ge 0\ge 0)}$).
\end{example}

\subsection{Determinant form}
One can easily provide explicit determinantal expressions for the 
coefficients $\gamma_{\ualpha,\ubeta}$. The following form was
stated in \S\ref{intro}:

\begin{theorem}\label{detform}
Let $d$ be such that $\alpha_{d+1}=0$. Then
$$\gamma_{\ualpha,\ubeta}=\sum \det
\left[
\begin{pmatrix}
\alpha_i-\ell_{i+1}^i-\cdots -\ell_d^i \\
\beta_j+(i-j)+\ell_i^1+\dots+\ell_i^{i-1}-\ell_{i+1}^i-\cdots -\ell_d^i 
\end{pmatrix}
\right]_{1\le i,j\le d}$$
where the summation is over the $\binom d2$ integers $\ell_i^k$, $1\le k
<i\le d$, subject to the conditions
$$0\le \ell_{k+1}^k + \cdots + \ell_d^k \le\alpha_{k+1}\quad.$$
\end{theorem}

\begin{proof}
This is obtained by computing explicitly the coefficients
$C_{\ualpha}(b_1, \dots,b_d)$ appearing in Theorem~\ref{complsym}, and
interpreting $\sum_{\sigma\in S_d}(-1)^\sigma
C_{\ualpha}(b_{\sigma(1)},\dots, b_{\sigma(d)})$ as a determinant.
\end{proof}

\begin{example}\label{detex}
For three-row diagrams, Theorem~\ref{detform} computes 
$\gamma_{\ualpha,\ubeta}$ as
$$\sum_{0\le \ell_2^1+\ell_3^1\le \alpha_2}\,\,
\sum_{0\le \ell_3^2\le \alpha_3} \det
\begin{pmatrix}
\binom{\alpha_1-\ell_2^1-\ell_3^1}{\beta_1-\ell_2^1-\ell_3^1} &
\binom{\alpha_1-\ell_2^1-\ell_3^1}{\beta_2-1-\ell_2^1-\ell_3^1} &
\binom{\alpha_1-\ell_2^1-\ell_3^1}{\beta_3-2-\ell_2^1-\ell_3^1} \\
\binom{\alpha_2-\ell_3^2}{\beta_1+1+\ell_2^1-\ell_3^2} &
\binom{\alpha_2-\ell_3^2}{\beta_2+\ell_2^1-\ell_3^2} &
\binom{\alpha_2-\ell_3^2}{\beta_3-1+\ell_2^1-\ell_3^2} \\
\binom{\alpha_3}{\beta_1+2+\ell_3^1+\ell_3^2} &
\binom{\alpha_3}{\beta_2+1+\ell_3^1+\ell_3^2} &
\binom{\alpha_3}{\beta_3+\ell_3^1+\ell_3^2}
\end{pmatrix}$$

For example, the coefficient of the Schubert class corresponding to the 
diagram
\begin{center}
\includegraphics[scale=.6]{pictures/200}
\end{center}
in the Chern-Schwartz-MacPherson class of the Schubert cell of
\begin{center}
\includegraphics[scale=.6]{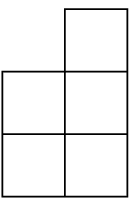}
\end{center}
is the sum of the determinants of the matrices{\tiny
\begin{align*}
&\begin{pmatrix}
1&0&0 \\
0&1&0 \\
0&1&1
\end{pmatrix}
\,,\,
\begin{pmatrix}
1&0&0 \\
0&0&0 \\
0&0&1 
\end{pmatrix}
\,,\,
\begin{pmatrix}
1&0&0 \\
0&1&0 \\
0&0&1
\end{pmatrix}
\,,\,
\begin{pmatrix}
1&0&0 \\
0&0&0 \\
0&0&0
\end{pmatrix}
\,,\,
\begin{pmatrix}
1&0&0 \\
0&1&0 \\
0&0&0
\end{pmatrix}
\,,\,
\begin{pmatrix}
1&0&0 \\
0&0&0 \\
0&0&0
\end{pmatrix}
\\
&\begin{pmatrix}
1&0&0 \\
0&2&1 \\
0&1&1
\end{pmatrix}
\,,\,
\begin{pmatrix}
1&0&0 \\
0&1&0 \\
0&0&1
\end{pmatrix}
\,,\,
\begin{pmatrix}
1&0&0 \\
0&2&1 \\
0&0&1
\end{pmatrix}
\,,\,
\begin{pmatrix}
1&0&0 \\
0&1&0 \\
0&0&0
\end{pmatrix}
\,,\,
\begin{pmatrix}
1&0&0 \\
0&1&2 \\
0&1&1
\end{pmatrix}
\,,\,
\begin{pmatrix}
1&0&0 \\
0&1&1 \\
0&0&1
\end{pmatrix}
\end{align*}}
that is,
$$1+0+1+0+0+0+1+1+2+0+(-1)+1=6\quad.$$
Duality implies that the same result may be obtained in a completely
different way: this coefficient must equal the coefficient of the
Schubert class corresponding to 
\begin{center}
\includegraphics[scale=.6]{pictures/11}
\end{center}
in the Chern-Schwartz-MacPherson class of the Schubert cell of
\begin{center}
\includegraphics[scale=.6]{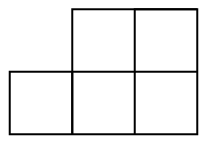}
\end{center}
that is, (again by Theorem~\ref{detform}) by the sum of binomial 
determinants
$$\sum_{\ell=0}^{\alpha_2}
\det
\begin{pmatrix}
\binom{\alpha_1-\ell}{\beta_1-\ell} & 
\binom{\alpha_1-\ell}{\beta_2-1-\ell} \\
\binom{\alpha_2}{\beta_1+1+\ell} & 
\binom{\alpha_2}{\beta_2+\ell}
\end{pmatrix}
=\det\begin{pmatrix}
3 & 1 \\ 1 & 2
\end{pmatrix}+
\det\begin{pmatrix}
1 & 0 \\ 0 & 1
\end{pmatrix}+
\det\begin{pmatrix}
0 & 0 \\ 0 & 0
\end{pmatrix}
=6\quad.$$
In general, duality implies rather complicated combinatorial identities 
between sums of binomial determinants. It would be interesting to give
a direct proof of these identities.
\end{example}

\subsection{Generating function}
It is also possible to present the numbers $\gamma_{\ualpha,\beta}$
directly as coefficients of the series expansion of a rational function; 
this is the form given in Theorem~\ref{main}.

\begin{theorem}\label{main2}
Let $d$ be such that $\alpha_{d+1}=0$. Then $\gamma_{\ualpha,\ubeta}$ 
equals the coefficient of 
$$t_1^{\alpha_1}\cdots t_d^{\alpha_d}\cdot 
u_1^{\beta_1}\cdots u_d^{\beta_d}$$
in the expansion of the rational function
$$\Gamma_d(\underline t,\underline u)=
\frac 1{(t_1^d\cdots t_d^1)(u_1^d\cdots u_d^1)}\cdot
\prod_{1\le i<j\le d} \frac{(t_i-t_j)(u_i-u_j)}{1-2t_j+t_i t_j}\cdot
\prod_{1\le i, j\le d}
\frac{1-t_i}{1-t_i(1+u_j)}$$
as a Laurent polynomial in $\Zbb[[t,u]]$.
\end{theorem}

This follows from a residue computation, based on the expression
given in Theorem~\ref{ratf}.

\begin{proof}
For the sake of notation we will give the argument for diagrams with
at most $d=2$ rows, that is,
$$\ualpha=(\alpha_1\ge \alpha_2\ge 0\ge \cdots)\quad;$$
the adaptations for larger $d$ present no difficulties.

For $d=2$, the expression in Theorem~\ref{ratf} is
$$\frac{(1+x_1)^{\alpha_1+1} (1+x_2)^{\alpha_2}}{1+x_1-x_2}\quad;$$
the number $C'_{\ualpha}(b_1,b_2)$ is the coefficient of 
$$x_1^{\alpha_1+2-b_1} x_2^{\alpha_2+1-b_2}\quad.$$
We can consider all $\ualpha$ at once as follows: 
$C'_{\ualpha}(b_1,b_2)$ is the coefficient of $x_1^{-b_1} x_2^{-b_2}$
in the Laurent expansion of
$$\frac{x_1^{-\alpha_1-2} x_2^{-\alpha_2-1}(1+x_1)^{\alpha_1+1} 
(1+x_2)^{\alpha_2}}{1+x_1-x_2}
=\frac{x_1^{-2}x_2^{-1}(1+x_1)}{1+x_1-x_2}
(1+x_1^{-1})^{\alpha_1}(1+x_2^{-1})^{\alpha_2}\quad;$$
that is, in the coefficient of $t_1^{\alpha_1}t_2^{\alpha_2}$
in the expansion of
$$\frac{x_1^{-2}x_2^{-1}(1+x_1)}{1+x_1-x_2}
\cdot \frac 1{1-t_1(1+x_1^{-1})}
\cdot \frac 1{1-t_2(1+x_2^{-1})}\quad.$$
The coefficient of $x_1^{-b_1}$ in this expression may be viewed as 
the residue 
$$\frac{x_2^{-1}}{1-t_2(1+x_2^{-1})}\,
\frac 1{2\pi i}
\oint
\frac{x_1^{b_1-2}(1+x_1)}{1+x_1-x_2}
\cdot \frac 1{1-t_1(1+x_1^{-1})}
\cdot \frac {dx_1}{x_1}\quad,$$
or, after the change of variable $y=(1-t_1(1+x_1^{-1}))x_1$, 
$$\frac{x_2^{-1}}{1-t_2(1+x_2^{-1})}\,
\frac 1{2\pi i}
\oint
\left(\frac{y+t_1}{1-t_1}\right)^{b_1-2}
\frac{1+y}{1+y-x_2(1-t_1)}\cdot
\frac {dy}{y(1-t_1)}\quad.$$
If $b_1\ge 2$, this is evaluated as
\begin{equation*}
\tag{*}
\frac{x_2^{-1}}{1-t_2(1+x_2^{-1})}\cdot
\frac{t_1^{b_1-2}}{(1-t_1)^{b_1-1}}\cdot
\frac 1{1-x_2(1-t_1)}\quad;
\end{equation*}
if $b_1=1$, an extra residue is picked up at $y=-t_1$, and is
evaluated as
$$\frac{x_2^{-1}}{1-t_2(1+x_2^{-1})}\cdot
\frac{-1}{t_1(1-x_2)}\quad;$$
but we can ignore this term, since we are only interested in
the coefficients of $t_1^{\alpha_1}$ for $\alpha_1\ge 0$.

An entirely analogous computation evaluates the coefficient of
$x_2^{-b_2}$ in (*) as
$$\frac{t_1^{b_1-2}}{(1-t_1)^{b_1-1}}\cdot
\frac{t_2^{b_2-1}}{(1-t_2)^{b_2-1}}\cdot
\frac 1{1-2t_2+t_1 t_2}\quad.$$

Summarizing, $C'_{\ualpha}(b_1,b_2)$ equals the coefficient of
$t_1^{\alpha_1}t_2^{\alpha_2}$ in the Laurent expansion of
$$\frac 1{1-2t_2+t_1 t_2}\cdot
\frac{(1-t_1)(1-t_2)}{t_1^2 t_2}\cdot
\frac{t_1^{b_1}}{(1-t_1)^{b_1}}\cdot
\frac{t_2^{b_2}}{(1-t_2)^{b_2}}\quad.$$
Multiplying by bookkeeping terms $u_1^{b_1}$, $u_2^{b_2}$
and adding over all $b_1\ge 0$, $b_2\ge 0$ shows that 
$C'_\ualpha(b_1,b_2)$ is the coefficient of 
$t_1^{\alpha_1}\,t_2^{\alpha_2}\,u_1^{b_1}\,u_2^{b_2}$
in the expansion of
$$\frac 1{1-2t_2+t_1 t_2}\cdot
\frac{(1-t_1)^2(1-t_2)^2}{t_1^2 t_2}\cdot
\frac 1{1-t_1(1+u_1)}\cdot
\frac 1{1-t_2(1+u_2)}\quad,$$
for all positive $b_1$, $b_2$. It follows that 
$\gamma_{\ualpha,\ubeta}$ equals the coefficient of 
$t_1^{\alpha_1}\,t_2^{\alpha_2}\,u_1^{b_1}\,u_2^{b_2}$ in
the expansion of
$$\frac 1{1-2t_2+t_1 t_2}\cdot
\frac{(1-t_1)^2(1-t_2)^2}{t_1^2 t_2}\cdot
\det\left[
\frac 1{1-t_i(1+u_j)}
\right]_{1\le i,j\le 2}$$
for $b_1=\beta_1+2$, $b_2=\beta_2+1$. That is, the coefficient of
$t_1^{\alpha_1}\,t_2^{\alpha_2}\,u_1^{\beta_1}\,u_2^{\beta_2}$
in the expansion~of
$$\frac 1{1-2t_2+t_1 t_2}\cdot
\frac{(1-t_1)^2(1-t_2)^2}{(t_1^2 t_2) (u_1^2 u_2)}\cdot
\det\left[
\frac 1{1-t_i(1+u_j)}
\right]_{1\le i,j\le 2}\quad.$$
The determinant may be evaluated directly; in general,
$$\det \left[
\frac 1{1-t_i(1+u_j)}
\right]_{1\le i,j\le d}
=\frac{\prod_{1\le i<j\le d}(t_i-t_j)(u_i-u_j)}
{\prod_{1\le i,j\le d} (1-t_i(1+u_j))}\quad,$$
an application of {\em Cauchy's double alternant\/}
(see for example \cite{MR1701596}, (2.7)).
This yields the formula given in the statement.
\end{proof}

\subsection{Contribution of one-row diagrams}
As an application of Theorem~\ref{main2}, we can give a rather
compact evaluation of the contribution of Schubert classes corresponding
to one-row diagrams to the CSM class of an arbitrary $\ualpha$; this
formula was mentioned in \S\ref{intro}.

\begin{corol}\label{onerow}
For all partitions $\ualpha$,
$$\sum_{r\ge 0} \gamma_{\ualpha,(r)}\, u^r
=\prod_{i\ge 1} (1+iu)^{\alpha_i-\alpha_{i+1}}\quad.$$
\end{corol}

\begin{proof}
Assume $\alpha_{d+1}=0$. By Theorem~\ref{main2}, $\gamma_{\ualpha,(r)}$ 
equals the coefficient of $t_1^{\alpha_1}\cdots t_d^{\alpha_d}\cdot 
u_1^r$ in the expansion of $\Gamma_d(\underline t,\underline u)$. 
Standard manipulations show that this equals the coefficient of 
$$t_1^{\alpha_1+d-1}\cdots t_{d-1}^{\alpha_{d-1}+1} t_d^{\alpha_d}
\cdot u^r$$
in the expansion of
$$F(\underline t)=\prod_{1\le i<j\le d}
\frac{t_i-t_j}{1-2t_j+t_i t_j}\cdot
\prod_{1\le i\le d} \frac 1{1-t_i(1+u)}\quad.$$
Now set $a_i=\alpha_i+d-i$, and note that $a_1>a_2>\cdots>a_d\ge 0$.
We say that two rational functions of $t_1,\dots,t_d$ {\em agree on 
the good cone\/} if the coefficients of $t_1^{a_1}\cdots t_d^{a_d}$
agree whenever $a_1>\cdots >a_d\ge 0$. We proceed to find a simpler
rational function agreeing with $F$ on the good cone.

A partial fraction decomposition of $F(\underline t)$ as a function
of $t_1$ gives
$$F(\underline t)=
\frac{C_1}{1-t_1(1+u)}+\sum_{j=2}^d \frac{C_j}{1-2t_j+t_1 t_j}$$
with $C_j$ rational functions in $t_2,\dots,t_d, u$, expanding
to elements of $\Zbb[[t_2,\dots,t_d,u]]$. It follows that
$F(\underline t)$ agrees with 
$$\frac{C_1}{1-t_1(1+u)}$$
on the good cone, since the exponent of $t_1$ in every other summand
is no larger than the exponent of $t_j$, for some $j>1$. 
Evaluating $C_1$ shows that $F(\underline t)$ and
$$\frac 1{1-t_1(1+u)}\cdot
\prod_{2\le i<j\le d}
\frac{t_i-t_j}{1-2t_j+t_i t_j}\cdot
\prod_{2\le i\le d} \frac 1{1+u-t_i(1+2u)}$$
agree on the good cone. 

Applying the same procedure to $C_1$, viewed as a function of $t_2$, 
shows that $F(\underline t)$ agrees with
$$\frac 1{(1-t_1(1+u))(1-t_2(1+2u))}\cdot
\prod_{3\le i<j\le d}
\frac{t_i-t_j}{1-2t_j+t_i t_j}\cdot
\prod_{3\le i\le d} \frac 1{1+2u-t_i(1+3u)}$$
on the good cone. Repeating again, at the $d$-th stage one obtains
that $F(\underline t)$ and 
$$G(\underline t):=\prod_{i=1}^d \frac 1{1+(i-1)u-t_i(1+iu)}$$
agree on the good cone. 

Therefore, $\gamma_{\ualpha,(r)}$ equals the coefficient of
$t_1^{\alpha_1+d-1}\cdots t_d^{\alpha_d}\cdot u^r$ in the expansion
of $G(\underline t)$. This can be evaluated easily, and yields
the statement.
\end{proof}

The formula given in Corollary~\ref{onerow} implies the following
pretty expression for the contribution of $[\Sbb((r))]$ to the
total Chern class of a Grassmannian $G_d(V)$:
this equals the coefficient of $u^r$ in
$$\frac 1{d!\, u^d} \sum_{i=0}^d \binom di (-1)^{d-i} 
(1+iu)^{\dim V}
\quad.$$
For example, the coefficients of one-row Schubert classes in
the total Chern class of $G_2(\Cbb^5)$ may be read off
$$\frac 1{2 u^2}\left(1-2(1+u)^5+(1+2u)^5\right)=10+30 u+35u^2+15u^3
\quad,$$
cf.~Example~\ref{g25}.


\section{Positivity and related issues}\label{pos}

\subsection{}\label{statcon}
Substantial computer experimentation suggests the following:

\begin{conj}\label{poscon}
For all $\ualpha$, $\csm(\Sbb(\ualpha)^\circ)\in 
A_*(\Sbb(\ualpha))$ is represented by an effective cycle.
\end{conj}

That is, we conjecture that the coefficients 
$\gamma_{\ualpha,\ubeta}$ are nonnegative, for all $\ubeta\le 
\ualpha$. 

As a consequence, the Chern-Schwartz-MacPherson class of every
Schubert {\em variety\/} would be positive in the sense of
\cite{MR85k:14004}, Chapter~12. This is well-known for
the Grassmannian (that is, for $\ualpha=(N^d)$): indeed, the
tangent bundle of a Grassmannian is generated by its sections
(\cite{MR85k:14004}, Examples~12.2.1 and~12.1.7). Such standard
positivity arguments do not seem to apply directly to 
Chern classes of singular varieties. Conjecture~\ref{poscon}, 
if true, may reflect a new positivity principle for such 
Chern classes.

Proving positivity from the explicit formulas obtained in
\S\ref{CSM} is an interesting challenge. Note that some of
the terms in the expression given in Theorem~\ref{complsym} 
may give a negative contribution to a coefficient 
(cf.~Example~\ref{exa22}: $x_1^2$~contributes $-e_{2,0}=-1$ 
to $\gamma_{(2\ge 2), (1\ge 1)}$); and some of the individual 
determinants in Theorem~\ref{detform} may be negative 
(cf.~Example~\ref{detex}). Conjecture~\ref{poscon} amounts to 
the statement that these negative contributions are 
canceled by corresponding positive contributions, and proving
that such cancellations always occur appears to be difficult
in general.

On the other hand, Corollary~\ref{onerow} implies that all
coefficients $\gamma_{\ualpha,(r)}$ of one-row Schubert classes in the
CSM class of any Schubert cell are positive; this is some evidence for
Conjecture~\ref{poscon}. In particular, Conjecture~\ref{poscon} is
immediate for Schubert cells corresponding to one-row diagrams. It is
also true for diagrams with $2$ rows: in this section we give two
different proofs of this easy result, which highlight two different
interesting facts on Chern-Schwartz-MacPherson classes of Schubert
cells.

These facts are a partial recursion for the coefficients 
$\gamma_{\ualpha,\ubeta}$, and a concrete interpretation of these 
coefficients as counting certain lattice paths between sets of points 
in the plane. We are able to prove positivity for all diagrams with~$3$ 
rows, by refining this latter counting argument. However, the proof in
this case is substantially more technical, and we will report on it 
elsewhere.

\subsection{Partial recursion}
Given a non-empty diagram $\ualpha$, we denote by $\ualpha^-$ the
diagram obtained by removing the rightmost column from $\ualpha$:
\begin{center}
\includegraphics[scale=.6]{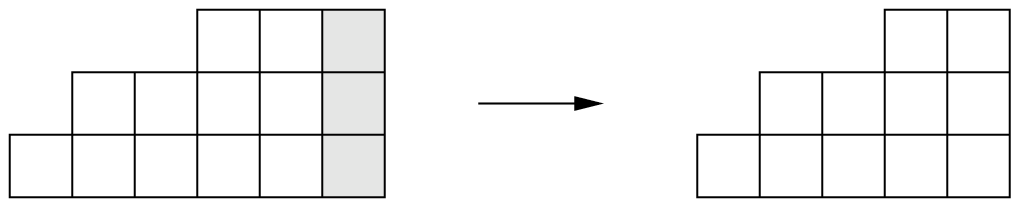}
\end{center}
That is, $\ualpha^-=\ualpha-(1^d)$, where $d$ is the largest
integer for which $\alpha_d\ne 0$.

We have the inclusion of Schubert {\em varieties\/}
$$\Sbb(\ualpha^-) \subset \Sbb(\ualpha)\quad.$$

\begin{prop}\label{adjform}
Let $\ualpha$ be a diagram with exactly $d$ rows, and embed
$\Sbb(\ualpha)$ in a Grassmannian $G_d(V)$. Let $\cS$ be the
universal subbundle on $G_d(V)$. Then
$$c_d(\cS^\vee)\cap \csm(\Sbb(\ualpha))
=c(\cS^\vee) \cap \csm(\Sbb(\ualpha^-))\quad,$$
and the same relation holds for CSM classes of cells.
\end{prop}

This statement may be interpreted as an `adjunction formula',
in the sense that the given identity is the formula that would
follow from the exact sequence of the normal bundle if
$\Sbb(\ualpha^-)\subset \Sbb(\ualpha)$ were nonsingular, 
with normal bundle the restriction of $\cS^\vee$. Simple
examples show that such formulas do not hold in general
for CSM classes; the fact that that they do hold
for the embeddings $\Sbb(\ualpha^-)\subset \Sbb(\ualpha)$
calls for a more geometric explanation.

\begin{proof}
By additivity, it suffices to prove the corresponding formula at
the level of CSM classes of Schubert cells:
$$c_d(\cS^\vee)\cap \csm(\Sbb(\ualpha)^\circ)
=c(\cS^\vee) \cap \csm(\Sbb(\ualpha^-)^\circ)\quad.$$
Let $\pi: \Vbb(\ualpha) \to \Sbb(\ualpha)$, $\pi^-: \Vbb(\ualpha^-)
\to \Sbb(\ualpha^-)$ be the birational isomorphisms constructed
in \S\ref{mainmapsec}. Recall that $\Vbb(\ualpha^-)$ embeds in 
$\Vbb(\ualpha)$, compatibly with these maps; and $[\Vbb(\ualpha^-)]$
has class $\xi_1\cdots \xi_d$ in $A_*(\Vbb(\ualpha))$ by 
Lemma~\ref{geomcla}.

By the very constructions of $\pi$,
$$\cS_1=\pi^* \cS\quad,$$
so that
$$c(\pi^* \cS)=\prod_{i=1}^d c(\cL_i^\vee)$$
(cf.~\S\ref{pb}). By Corollary~\ref{urform} and the projection 
formula,
\begin{align*}
c_d(\cS^\vee)\cap \csm(\Sbb(\ualpha)^\circ)
&=\pi_*\left(\left(\prod_{i=1}^d \xi_i\right) \prod_{i\ge 1}  
c(\cL_i^\vee)^{\alpha_i-\alpha_{i+1}}\, c(\cL_i^\vee \otimes 
\cQ_{i+1})\cap [\Vbb(\ualpha)]\right)\\
&=\pi_*\left(\prod_{i\ge 1} c(\cL_i^\vee)^{\alpha_i-
\alpha_{i+1}}\, c(\cL_i^\vee \otimes \cQ_{i+1})\cap [\Vbb(\ualpha^-)]
\right)
\end{align*}
Here $\cL_i$, $\cQ_i$ are the structure bundles on $\Vbb(\ualpha)$.
Denoting by $\cL_i^-$, $\cQ_i^-$ the structure bundles on 
$\Vbb(\ualpha^-)$ we have, as observed in \S\ref{strmaps}, 
$$\cL_i|_{\Vbb(\ualpha^-)}={\cL_i^-}^\vee
\quad,\quad
\cQ_i|_{\Vbb(\ualpha^-)}={\cQ_i^-}\oplus \cO\quad.$$
Therefore,
$$(\cL_i^\vee \otimes \cQ_{i+1})|_{\Vbb(\ualpha^-)}
=({\cL_i^-}^\vee \otimes \cQ_{i+1}^-)\oplus {\cL_i^-}^\vee\quad,$$
giving
\begin{align*}
c_d(\cS^\vee)\cap &\csm(\Sbb(\ualpha)^\circ)
=\pi^-_*\left(\prod_{i\ge 1} c({\cL_i^-}^\vee)^{\alpha_i-
\alpha_{i+1}}\, c({\cL^-_i}^\vee \otimes \cQ^-_{i+1})
\, c({\cL_i^-}^\vee)\cap [\Vbb(\ualpha^-)]\right)\\
&=\pi^-_*\left(\left(\prod_{i\ge 1} c({\cL_i^-}^\vee)\right) 
\prod_{i\ge 1} c({\cL_i^-}^\vee)^{\alpha_i-\alpha_{i+1}}\, 
c({\cL^-_i}^\vee \otimes \cQ^-_{i+1})\cap [\Vbb(\ualpha^-)]\right)\\
&=c(\cS^\vee)\cap \csm(\Sbb(\ualpha^-)^\circ)
\end{align*}
again by the projection formula and Corollary~\ref{urform}.
\end{proof}

By duality, Proposition~\ref{adjform} implies the following
alternative formulation. For a non-empty diagram $\ualpha$, 
denote by $\ualpha'$ the diagram obtained by removing the {\em
bottom row\/} from $\ualpha$:
\begin{center}
\includegraphics[scale=.6]{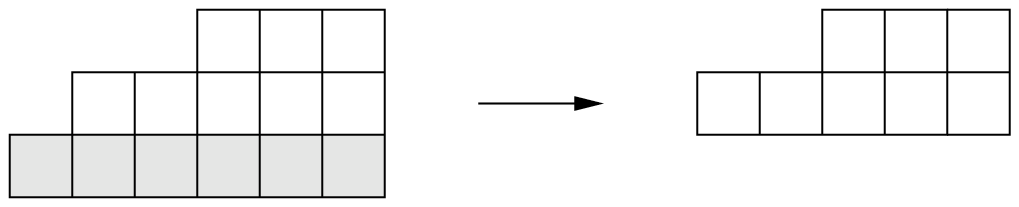}
\end{center}
That is, $\ualpha'$ is the truncation often used elsewhere in this
article, for inductive arguments.

\begin{prop}\label{adjform2}
Let $\ualpha$ be a nonempty diagram; embed $\Sbb(\ualpha)$ in the
Grassmannian $G_d(V)$, where $\dim V=\alpha_1+d$, and let $\cQ$ 
be the universal quotient bundle on $G_d(V)$. 
Then
$$c_{\alpha_1}(\cQ)\cap \csm(\Sbb(\ualpha))
=c(\cQ) \cap \csm(\Sbb(\ualpha'))\quad,$$
and the same relation holds for CSM classes of cells.
\end{prop}

\begin{proof}
Apply Proposition~\ref{adjform} to the transpose of $\ualpha$;
this gives the statement, since the isomorphism 
$G_d(V)\cong G_{\alpha_1}(V)$ interchanges
the universal quotient bundle with the (dual) universal subbundle.
\end{proof}

These formulas imply that {\em some\/} of the coefficients in the 
CSM class of a Schubert cell (or variety) are determined by the CSM 
class of Schubert cells for smaller diagrams. 
For example: provided the class $\csm(\Sbb(\ualpha^-)^\circ)$ is 
known, Proposition~\ref{adjform} determines the coefficients of 
the classes which survive the product by $c_d(\cS^\vee)$, where 
$d$ is the number of rows of the diagram $\ualpha$. 
These are the classes $[\Sbb(\ubeta)]$ for subdiagrams 
$\ubeta\le \ualpha$ also consisting of $d$ rows.

\begin{example}\label{331}
Let $\ualpha=(3\ge 3\ge 1)$:
\begin{center}
\includegraphics[scale=.6]{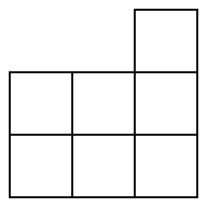}
\end{center}
Then the coefficients in $\csm(\Sbb(\ualpha)^\circ)$ of the six 
classes corresponding to subdiagrams of $(3\ge 3\ge 1)$ with exactly 
three rows are determined by Proposition~\ref{adjform} from
$\csm(\Sbb(\ualpha^-)^\circ)=\csm(\Sbb(2\ge 2)^\circ)$, that is
\begin{center}
\includegraphics[scale=.4]{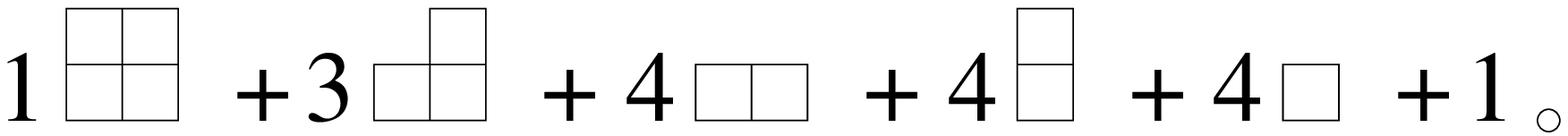}
\end{center}
(denoting each class by the corresponding diagram). 
Multiplying by the Chern class of the dual of the universal 
subbundle $\cS^\vee$ amounts to simple applications of Pieri's 
formula, and yields
\begin{center}
\includegraphics[scale=.4]{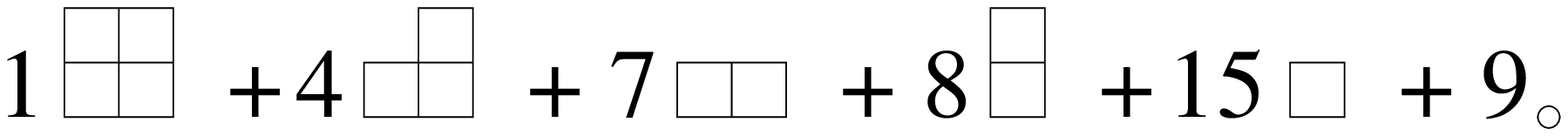}
\end{center}
According to Proposition~\ref{adjform}, the classes of diagrams
obtained by adding one vertical 3-column to each of these diagrams
appear with the same coefficients in $\csm(\Sbb(\ualpha)^\circ)$.
Similarly, applying Proposition~\ref{adjform2} allows us to 
recover the coefficients of classes corresponding to subdiagrams
with the same number of columns as $\ualpha$, provided one knows
$\csm(\Sbb(3\ge 1)^\circ)$, that is
\begin{center}
\includegraphics[scale=.4]{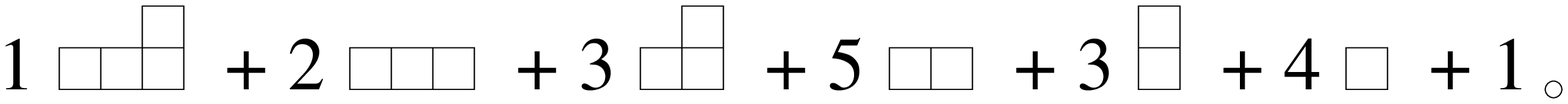}
\end{center}
Multiplying by $c(\cQ)$ gives
\begin{center}
\includegraphics[scale=.4]{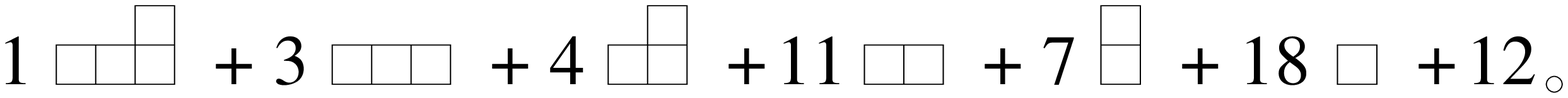}
\end{center}
giving the coefficients of classes for diagrams obtained from
these by adding one 3-row. The conclusion is that 
$\csm(\Sbb(\ualpha)^\circ)$ must equal
\begin{center}
\includegraphics[scale=.35]{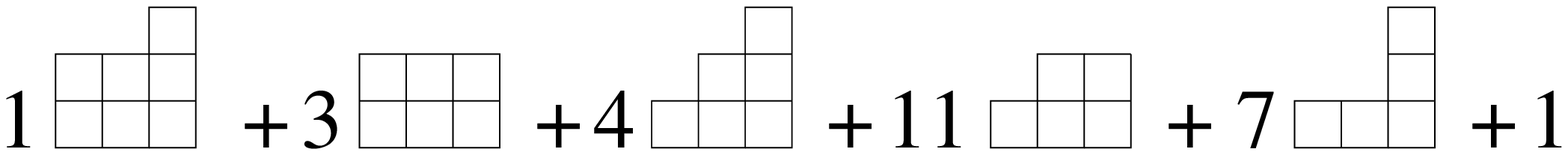}
\end{center}
plus classes for diagrams contained in the $2\times2$ rectangle.
That is, recursion suffices to compute most coefficients of
$\csm(\Sbb(\ualpha)^\circ)$.
\end{example}

\subsection{Positivity for two rows: first proof}
The relevance of these considerations to the positivity question
is the following: assuming that positivity has been established for 
all diagrams strictly smaller than $\ualpha$, then 
Propositions~\ref{adjform} and~\ref{adjform2} imply that the 
coefficients of $[\Sbb(\ubeta)]$ in $\csm(\Sbb(\ualpha)^\circ)$ 
are nonnegative, for all $\ubeta\le \ualpha$ with the same number
of rows or the same number of columns. Indeed, multiplication by
$c(\cS^\vee)$ or $c(\cQ)$ preserve effectivity. Therefore:

\begin{lemma}\label{recclaim}
In order to prove Conjecture~\ref{poscon}, it suffices to show that 
$\gamma_{\ualpha,\ubeta}\ge 0$ for all $\ualpha$ and all $\ubeta
< \ualpha$ with strictly fewer rows and columns than $\ualpha$.
\end{lemma}

The positivity of CSM classes for two-row diagrams follows
immediately from Lemma~\ref{recclaim} and Corollary~\ref{onerow}.
Indeed, positivity for one-row diagrams holds (as pointed out
in \S\ref{statcon}); by Lemma~\ref{recclaim}, verifying
the conjecture for two-row diagrams $\ualpha$ amounts to verifying 
that $\gamma_{\ualpha,\ubeta}\ge 0$ for all one-row diagrams
$\ubeta$, and Corollary~\ref{onerow} shows that this is indeed
the case: for $\ualpha=(\alpha_1\ge \alpha_2)$ and $\ubeta=(r)$,
$$\sum_{r\ge 0} \gamma_{\ualpha,(r)}\, u^r
=(1+u)^{\alpha_1-\alpha_2} (1+2u)^{\alpha_2}$$
has nonnegative coefficients.

\subsection{Nonintersecting lattice paths, and second proof
of positivity for two rows}
Positivity for cells corresponding to two-row diagrams is also
immediate from the following concrete interpretation of the 
coefficients $\gamma_{\ualpha,\ubeta}$ in this case: for 
$\ualpha=(\alpha_1\ge \alpha_2)$, $\ubeta=(\beta_1\ge\beta_2)$, 
with $\beta\le \alpha$, {\em $\gamma_{\ualpha,\ubeta}$ equals the 
number of certain non-intersecting lattice paths joining pairs 
of points in the plane.\/}

Here is the precise statement. A {\em lattice path\/} joining
two points $A$, $B$ of the integer lattice in the plane is a
sequence of lattice points starting from $A$ and ending in $B$,
such that each point is one step to the right or down from the
preceding one. Given $\alpha_1\ge \alpha_2$, consider the 
two sets of $\alpha_2+1$ points in the plane 
$$A_1=\{A_1^\ell\}_{0\le \ell \le \alpha_2} \quad,\quad
A_2=\{A_2^\ell\}_{0\le \ell \le \alpha_2}\quad,$$
where
$$A_1^\ell = (\ell+2,\alpha_1+2) \quad , \quad
A_2^\ell=(1-\ell,\alpha_2+1-\ell)\quad.$$
\begin{center}
\includegraphics[scale=.8]{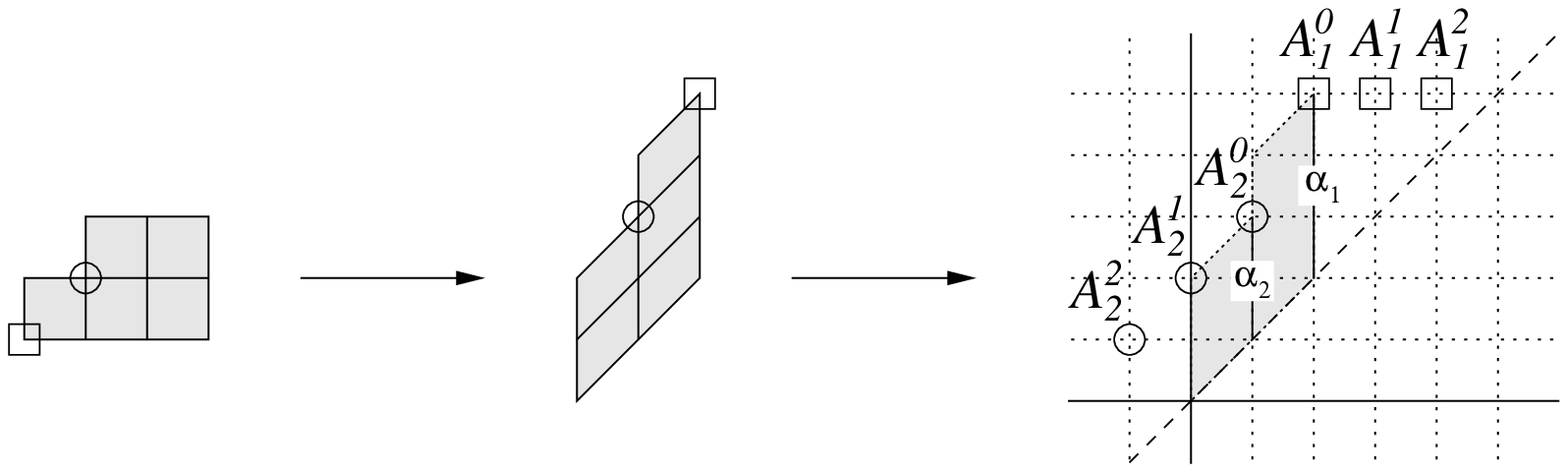}
\end{center}
Given $\beta_1\ge \beta_2$, also consider the pair of points
$$B_1=(\beta_1+2,\beta_1+2)\quad,\quad B_2=(\beta_2+1,\beta_2+1)$$
on the main diagonal. 
Let $\Pi(\ell)$, for $0 \leqslant \ell \leqslant \alpha_2$, be the
set of pairs $(\pi_1,\pi_2)$ such that $\pi_r$ is a path from 
$A_r^\ell$ to $B_r$ and $\pi_1$ doesn't intersect $\pi_2$.  

\begin{theorem}\label{latticepaths}
With notation as above, the coefficient $\gamma_{\ualpha,\ubeta}$ is
\[ \sum_{l=0}^{\alpha_2} |\Pi(\ell)| \] where $|\Pi(\ell)|$ denotes the 
cardinality of the set $\Pi(\ell)$.
\end{theorem}

\begin{example}
The contribution of $[\Sbb(1\ge 1)]$ to $\csm(\Sbb(3\ge 2)^\circ)$
is $\gamma_{3\ge 2,1\ge 1}=6$; the $6$ pairs of paths are
\begin{center}
\includegraphics[scale=.6]{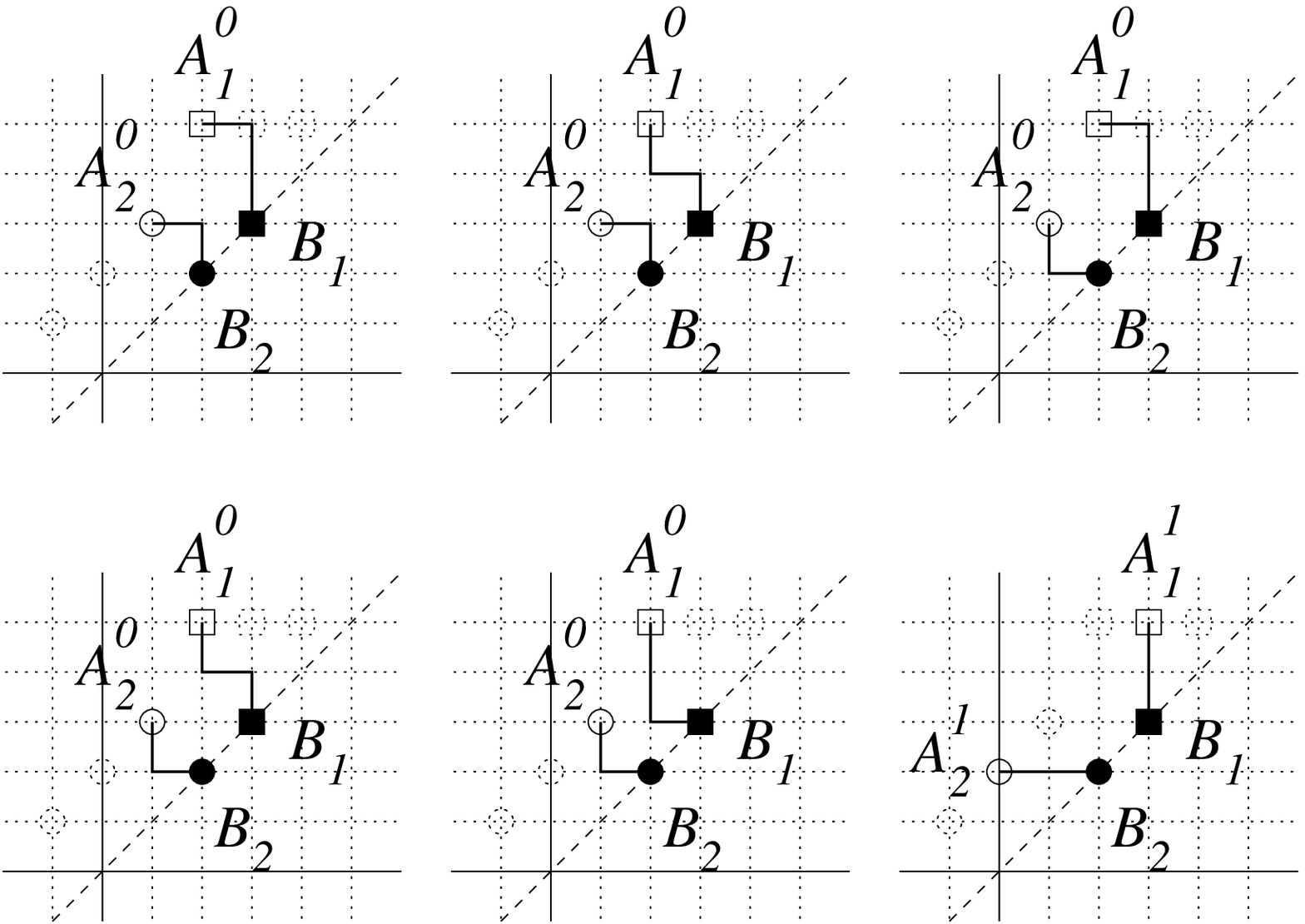}
\end{center}
There are 5 nonintersecting paths from the first pair of points 
in $A_1$, $A_2$ to $B_1$, $B_2$; 1 path from the second pair,
and no paths from the third pair. One pair of lattice paths
$A_1^0 \to B_1$, $A_2^0 \to B_2$:
\begin{center}
\includegraphics[scale=.5]{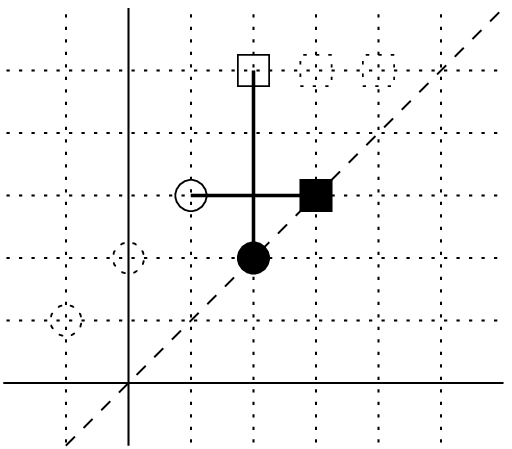}
\end{center}
is not counted since the paths meet.
\end{example}

\begin{proof}
Note that $A_2$ lies entirely due south-west of $A_1$, both are north 
of the main diagonal, and $B_2$ is south-west of $B_1$. This implies 
that every lattice path $A_1^\ell \to B_2$ meets every lattice path
from $A_2^\ell \to B_1$. Under this hypothesis, the number
of non-intersecting paths considered in the statement equals
$$\sum_{\ell=0}^{\alpha_2}
\det\left[
\begin{matrix}
P(A^\ell_1 \to B_1) & P(A^\ell_1 \to B_2) \\
P(A^\ell_2 \to B_1) & P(A^\ell_2 \to B_2)
\end{matrix}
\right]$$
where $P(A\to B)$ denotes the number of lattice paths from $A$ to
$B$; this follows from the `Lindstr\"om-Gessel-Viennot theorem'
(Theorem~6 in \cite{MR2178686}).

Now recall that, for $a,b$ nonnegative integers, the binomial coefficient 
$\binom ab$ equals the number of lattice paths joining the points
$(0,a)$ and $(b,b)$. Applying evident translations,
\begin{alignat*}{3}
P(A^\ell_1 \to B_1) &=\binom{\alpha_1-\ell}{\beta_1-\ell}\quad &,\quad
& P(A^\ell_1 \to B_2)=\binom{\alpha_1-\ell}{\beta_2-1-\ell} \\
P(A^\ell_2 \to B_1) &=\binom{\alpha_2}{\beta_1+1+\ell}\quad &,\quad
& P(A^\ell_2 \to B_2)=\binom{\alpha_2}{\beta_2+\ell}
\end{alignat*}
Therefore, the number of nonintersecting paths given in the statement
is
$$\sum_{\ell=0}^{\alpha_2}
\det\left[
\begin{matrix}
\binom{\alpha_1-\ell}{\beta_1-\ell} & 
\binom{\alpha_1-\ell}{\beta_2-1-\ell} \\
\binom{\alpha_2}{\beta_1+1+\ell} & 
\binom{\alpha_2}{\beta_2+\ell}
\end{matrix}
\right]\quad.$$
By Theorem~\ref{detform} this equals $\gamma_{\ualpha,\ubeta}$,
concluding the proof.
\end{proof}




\begin{thebibliography}{Sch65b}

\bibitem[Alu]{math.AG/0507029}
P.~Aluffi.
\newblock {Limits of Chow groups, and a new construction of
  Chern-Schwartz-MacPherson classes}.
\newblock FSU05-14, arXiv:math.AG/0507029.

\bibitem[Alu99]{MR1717120}
Paolo Aluffi.
\newblock Differential forms with logarithmic poles and
  {C}hern-{S}chwartz-{M}ac{P}herson classes of singular varieties.
\newblock {\em C. R. Acad. Sci. Paris S\'er. I Math.}, 329(7):619--624, 1999.

\bibitem[Alu06]{MR2209219}
Paolo Aluffi.
\newblock Classes de {C}hern des vari\'et\'es singuli\`eres, revisit\'ees.
\newblock {\em C. R. Math. Acad. Sci. Paris}, 342(6):405--410, 2006.

\bibitem[Dem74]{MR0354697}
Michel Demazure.
\newblock D\'esingularisation des vari\'et\'es de {S}chubert
  g\'en\'eralis\'ees.
\newblock {\em Ann. Sci. \'Ecole Norm. Sup. (4)}, 7:53--88, 1974.
\newblock Collection of articles dedicated to Henri Cartan on the occasion of
  his 70th birthday, I.

\bibitem[Ful84]{MR85k:14004}
William Fulton.
\newblock {\em Intersection theory}.
\newblock Springer-Verlag, Berlin, 1984.

\bibitem[GP02]{MR1893006}
Mark Goresky and William Pardon.
\newblock Chern classes of automorphic vector bundles.
\newblock {\em Invent. Math.}, 147(3):561--612, 2002.

\bibitem[Kra99]{MR1701596}
C.~Krattenthaler.
\newblock Advanced determinant calculus.
\newblock {\em S\'em. Lothar. Combin.}, 42:Art. B42q, 67 pp. (electronic),
  1999.
\newblock The Andrews Festschrift (Maratea, 1998).

\bibitem[Kra05]{MR2178686}
C.~Krattenthaler.
\newblock Advanced determinant calculus: a complement.
\newblock {\em Linear Algebra Appl.}, 411:68--166, 2005.

\bibitem[Mac74]{MR50:13587}
Robert~D. MacPherson.
\newblock Chern classes for singular algebraic varieties.
\newblock {\em Ann. of Math. (2)}, 100:423--432, 1974.

\bibitem[PP95]{MR1311826}
Adam Parusi{\'n}ski and Piotr Pragacz.
\newblock Chern-{S}chwartz-{M}ac{P}herson classes and the {E}uler
  characteristic of degeneracy loci and special divisors.
\newblock {\em J. Amer. Math. Soc.}, 8(4):793--817, 1995.

\bibitem[Sch65a]{MR35:3707}
M.-H. Schwartz.
\newblock Classes caract\'eristiques d\'efinies par une stratification d'une
  vari\'et\'e analytique complexe. {I}.
\newblock {\em C. R. Acad. Sci. Paris}, 260:3262--3264, 1965.

\bibitem[Sch65b]{MR32:1727}
M.-H. Schwartz.
\newblock Classes caract\'eristiques d\'efinies par une stratification d'une
  vari\'et\'e analytique complexe. {II}.
\newblock {\em C. R. Acad. Sci. Paris}, 260:3535--3537, 1965.

\bibitem[Vak]{math.AG/0302294}
Ravi Vakil.
\newblock {A geometric Littlewood-Richardson rule}, arXiv:math.AG/0302294.
\newblock to appear in the Annals of Math.

\end{thebibliography}
\end{document}